\documentclass[12pt]{article}

\usepackage{latexsym, longtable}
\usepackage{amssymb, amsmath, xspace, lscape,  latexsym}

\renewcommand{\dim}{\mathrm{dim \,}}

\renewcommand{\le}{\leqslant}
\renewcommand{\ge}{\geqslant}
\renewcommand{\natural}{\sharp}

\newcommand{\ve}{\varepsilon}

\def\beq#1#2\eeq{%
        \begin{equation}%
        \label{#1}%
            #2%
        \end{equation}%
    }

\newcommand{\mref}[1]{(\ref{#1})}

\newcommand{\p}{\partial}

\renewcommand{\a}{\alpha}
\renewcommand{\b}{\beta}
\newcommand{\g}{\gamma}

\newcommand{\C}{\mathbb C}
\newcommand{\N}{\mathbb N}
\newcommand{\Z}{\mathbb Z}

\renewcommand{\hat}{\widehat}
\renewcommand{\tilde}{\widetilde}

\def\btheor#1\etheor{%
        \begin{theor}%
            #1%
        \end{theor}%
    }

    \def\bsled#1\esled{%
        \begin{sled}%
            #1%
        \end{sled}%
    }

\newtheorem{theorem}{Theorem}

\newtheorem{lemma}{Lemma}
\newtheorem{prop}{Proposition}

\newtheorem{remark}{Remark}

\newcommand{\nad}[2]{\genfrac{}{}{0pt}{}{#1}{#2}}

\def\hm#1{#1\nobreak\discretionary{}{\hbox{\m@th$#1$}}{}}
\def\mi#1{\discretionary{\hbox{\m@th$#1$}}{\hbox{\m@th$#1$}}{}}

\vspace{4ex}
\begin{document}

\begin{center} {\Large \bf
Generalized Calogero-Moser systems from rational Cherednik algebras
}
\end{center}

\vspace{3mm}

\begin{center}
{\bf M.V. Feigin}
\end{center}

\vspace{3mm}

\noindent
School of Mathematics and Statistics, University of Glasgow,
15 University Gardens, Glasgow G12 8QW, UK. Email: misha.feigin@glasgow.ac.uk

\vspace{9mm}

\begin{abstract} We consider ideals of polynomials vanishing on the $W$-orbits
of the intersections of mirrors of a finite reflection group $W$.
We determine all such ideals which are invariant
under the action of the corresponding rational Cherednik algebra
hence form submodules in the
polynomial module. We show that a quantum integrable system can be
defined for every such ideal for a real reflection group $W$. This
leads to known and new integrable systems of Calogero-Moser type
which we explicitly specify. In the case of classical Coxeter groups we also obtain generalized Calogero-Moser systems with added quadratic potential.
\end{abstract}



\section{Introduction}

The usual Calogero-Moser (CM) system describes a pairwise
interaction of $n$ particles on the line with the inverse square
potential. This system appeared and was studied in \cite{Calogero1},
\cite{Calogero2}, \cite{Moser}. Olshanetsky and Perelomov introduced
CM systems related to the root systems of finite Coxeter groups
\cite{OP}. An elegant and uniform proof of integrability was
proposed by Heckman \cite{Heckman} who used Dunkl operators
\cite{Dunkl}. After that quantum CM systems related to the root
systems became ultimately related to the rational Cherednik algebras
$H_{1,c}$ \cite{EG}. More exactly, CM operators and their quantum
integrals can be thought of as elements of the spherical subalgebra
of the rational Cherednik algebra. We refer to the book
\cite{Etingof} for the exposition of this and other developments.

Generalized CM systems related to non-symmetric arrangements of
hyperplanes appeared in the work of Chalykh, Veselov and the author
\cite{CFV1}. These systems were studied in \cite{CFV2}, \cite{CFV3}
where integrability was established with the help of Baker-Akhiezer
functions. Two families of operators were found corresponding to the
deformations of the root systems $\mathcal A_n(m)$, $\mathcal
C_n(m,l)$. The $\mathcal A_n(m)$ system describes Calogero-Moser
type pairwise interaction where one particle has different mass from
the mass of all other particles.

Sergeev considered the generalization of $\mathcal A_n(m)$ CM system
with arbitrary number of particles of each of the two types
\cite{Sergeev}. He obtained the corresponding
Calogero-Moser-Sutherland (CMS) operator at $m=-1/2$ as the radial
part of the Laplace-Beltrami operator on a symmetric superspace
\cite{Sergeev, Sergeev2}. Sergeev and Veselov introduced further
generalization of the families of CMS operators by defining special
deformations of the root systems of contragredient  superalgebras
Lie \cite{SV1}. Then integrability for the classical series was
established in \cite{SV1} by rather involved computations.

In the works \cite{SV2,SV3} Sergeev and Veselov gave another proof
of integrability of CMS systems related to deformations of classical
root systems. They showed that the generalized CMS operators are
restrictions to generalized discriminants of the non-deformed CMS
operators of $A$ or $BC$ type acting in the space of symmetric
functions in infinite number of variables. Special eigenfunctions of
these systems were studied by Halln\"as and Langmann in
\cite{HL,H3}.

In the present paper we approach integrability of generalized CM
systems making use of Dunkl operators and special representations of
rational Cherednik algebras. We show that for special values of
parameters the Dunkl operators can be restricted to certain
parabolic strata which are the Coxeter orbits of intersections of
the mirrors. Equivalently, the ideals of polynomials vanishing on
these strata are submodules in the polynomial representation of the
corresponding rational Cherednik algebra. Then sum of squares of base
Dunkl operators takes the form of a generalized CM operator when acting
on the invariants. In this way we recover CM systems from
\cite{CFV1}-\cite{SV3} at the special values of parameters and also
obtain new integrable generalized CM systems. It is simpler to work
in the non-symmetric settings and we give complete description of
the CM systems which can be obtained in this way. The corresponding
parabolic strata are easy to describe. Namely, in the case of
constant multiplicity the Coxeter subgraph for the parabolic
subgroup defining the stratum should have the same Coxeter number
for all its connected components.

We note that the interest in the submodules of the polynomial
representation of the rational Cherednik algebras goes back to the
pioneering work by Dunkl, de Jeu and Opdam \cite{DJO} where the set
of singular parameters when the representation is reducible was
completely determined. The singular set for the trigonometric
degeneration of the Cherednik algebra was found by Etingof in
\cite{Et2}, the non-degenerate case was studied by Cherednik in
\cite{Cherednik}. The actual submodules for the non-degenerate
Cherednik algebras were under investigation, in particular, by
Kasatani for $A$ and $C$ cases in \cite{K1}, \cite{K2},  and by
Cherednik in \cite{Cherednik}. The representations of the Coxeter
group on the singular vectors were determined by Dunkl in the
rational $A$ case \cite{dun}.

The paper is organized as follows. In Section \ref{sec1} we
determine parabolic strata of the finite real reflection groups such
that the ideals of polynomials vanishing on these strata are
submodules for the polynomial representation of the corresponding
rational Cherednik algebra.
These strata are rational Coxeter versions of Kasatani's
nonsymmetric vanishing conditions (\cite{K1}, see also \cite{BL0,
BL}). It follows from \cite{Gi}, \cite{BE} that it is necessary to
impose vanishing conditions on the parabolic strata in order to get
non-trivial invariant ideals.

In Section \ref{sec2} we show that parabolic strata defining
invariant ideals also determine quantum integrable systems of CM
type corresponding to the sets of vectors which are obtained by
projections of the original Coxeter root system. In Section
\ref{sec3} we explicitly specify these invariant parabolic strata
and the corresponding CM systems in the case of classical Coxeter
groups. The corresponding generalized CM systems are known to be
integrable by \cite{CFV3}, \cite{SV1}. In Section \ref{sec4} we show
that systems found in Section \ref{sec3} remain integrable if a
quadratic term is added to the Hamiltonian. For that we review the
proof of integrability of the CM systems for classical series in the
external quadratic field through the Dunkl operators \cite{Poly}.
Then our restriction procedure can be applied in this case as well.
In Section \ref{sec5} we explicitly determine generalized CM systems
corresponding to the invariant parabolic ideals for the exceptional
Coxeter groups. In Section~\ref{sec6} we determine invariant
parabolic ideals for complex reflection groups.

\section{Invariant parabolic ideals}\label{sec1}
\label{section-real-parabolic}

Let $W$ be a finite real reflection group acting by orthogonal
transformations in its complexified reflection representation
$V=\C^N$. Let $\mathcal R$ be the corresponding Coxeter root
system, and let $\Gamma$ be the corresponding Coxeter graph (see,
e.g., \cite{Humphreys}). We assume that a positive subsystem
$\mathcal R_+ \subset \mathcal R$ is chosen so the vertices of the
graph $\Gamma$ can be identified with the simple roots. Similarly,
for a subgraph $\Gamma_0 \subset \Gamma$ we will denote by
$\Gamma_0^v$ the set of roots corresponding to the vertices of
$\Gamma_0$.

Let $c(\a)=c_\a$ be a $W$-invariant function on the set of roots
$\mathcal R$. The {\it rational Cherednik algebra}
$H_c=H_c^{\mathcal R}$ is associated with the root system $\mathcal
R$ and multiplicity $c$ (see \cite{EG}; in this paper we assume that
parameter $t=1$). Also in this paper we will need only the faithful
representation of $H_c$ in the space of polynomials
$\C[x]=\C[x_1,\ldots,x_N]$. Any element of $H_c$ acts on $p\in\C[x]$
as a linear combination of the compositions $r(\nabla) w q (x)p$ where $q(x) \in \C[x]$, $w\in
W$, and $r(\nabla)=r(\nabla_1,\ldots,\nabla_N)$ is a polynomial in
Dunkl operators corresponding to the bases directions
$e_1,\ldots,e_N$. For any direction $\xi \in \C^N$ the Dunkl
operator $\nabla_\xi$ is defined as \beq{dunkldef}  \nabla_\xi =
\p_\xi - \sum_{\a\in {\mathcal R}_+} \frac{c_\a
(\a,\xi)}{(\a,x)}(1-s_\a), \eeq where $(\cdot,\cdot)$ is the
standard scalar product in $V$, and $s_\a$ is the orthogonal
reflection with respect to the hyperplane $(\a,x)=0$. Note that the
Dunkl operators satisfy commutativity $[\nabla_\xi,\nabla_\eta]=0$
(\cite{Dunkl}) and it is clear that $\nabla_\xi \C[x] \subset
\C[x]$.

Let $\Gamma_0$ be a subgraph of the Coxeter graph $\Gamma$ obtained
by specifying some vertices of $\Gamma$ and preserving all the edges
between these vertices. The subgraph $\Gamma_0$ defines the plane
$\pi$ obtained as the intersection of the corresponding mirrors
\beq{pidef} \pi = \{x\in V \mid (\b,x)=0 \quad \forall \b\in \Gamma_0^v
\}. \eeq The associated {\it parabolic stratum} is defined as
$$
D_{\Gamma_0}=\bigcup_{w \in W} w(\pi).
$$
We define the corresponding {\it parabolic ideal} $I_{\Gamma_0}$ as
a set of polynomials vanishing on the stratum, that is
$I_{\Gamma_0}=\{p\in \C[x]\mid p|_{D_{\Gamma_0}}=0\}$. It is obvious
that $I_{\Gamma_0}$ is an ideal in $\C[x]$ and that it is
$W$-invariant. We are going to determine the parabolic strata
$D_{\Gamma_0}$ which define ideals $I_{\Gamma_0}$ invariant under
the whole rational Cherednik algebra $H_c$.

\begin{theorem}\label{mth1}
Let $\Gamma_0=\coprod_{i=1}^k \Gamma_i$ be the decomposition of the
subgraph $\Gamma_0 \subset \Gamma$ into the connected components.
Then the parabolic ideal $I_{\Gamma_0}$ is invariant under the
algebra $H_c$ if and only if for any $i=1,\ldots,k$ we have
\beq{idealid} \sum_{\a \in V_i\bigcap \mathcal R} \frac{
c_\a(\a,u)(\a,v)}{(\a,\a)}= (u,v) \eeq for any $u,v \in V_i$, where
$V_i$ is a linear space spanned by the roots $\Gamma_i^v$.
\end{theorem}
{\bf Proof.} Denote by $V_0= \bigoplus_{i=1}^k V_i$. Let $f \in
I_{\Gamma_0}$. We are going to analyze the submodule condition
$\nabla_\xi f|_{D_{\Gamma_0}}=0$ where $ \nabla_\xi$ is the Dunkl
operator \mref{dunkldef}. At first we consider the condition
$\nabla_\xi f|_{\pi}=0$. For that we recall that $f|_\pi=0$ and we
can represent polynomial $f$ in the form
$$
f= \sum_{\b \in \Gamma_0^v} f_\b(x) (\b,x),
$$
where $f_\b$ are some polynomials. Since $s_\a f|_\pi =
f|_{s_\a\pi}=0$ we note that $\frac{1-s_\a}{(\a,x)}f|_\pi=0$ if
$(\a,x)|_\pi \ne 0$. Therefore we rearrange
$$
\nabla_\xi f |_\pi = \sum_{\b \in \Gamma_0^v} \left((\b,\xi) f_\b
|_\pi  - \sum_{\a\in  {\mathcal R}_+\bigcap V_0} \frac{2 c_\a (\a,
\xi) (\a,\b)}{(\a,\a)} f_\b|_\pi\right).
$$
By collecting the coefficients at $f_\b$ it follows that
$$
(\b,\xi)- \sum_{\a \in  {\mathcal R}_+\bigcap V_0} \frac{2 c_\a (\a,
\xi) (\a,\b)}{(\a,\a)} =0,
$$
which is equivalent to the property \mref{idealid} as $V_i$ are
pairwise orthogonal for $i>0$.

Conditions $\nabla_\xi f |_{w(\pi)}=0$ for non-trivial $w \in W$ are
obtained from conditions \mref{idealid} by the $W$-action, they are
equivalent to the properties \mref{idealid} hence the theorem is
proven.

\begin{theorem}\label{mth2}
Assume that the multiplicity function $c$ is constant on the roots
$\Gamma_0^v$ of the subgraph $\Gamma_0 \subset \Gamma$. Then the
ideal $I_{\Gamma_0}$ is $H_c$-invariant if and only if all the
connected components of $\Gamma_0$ have same Coxeter number $h=1/c$.
\end{theorem}

This theorem is a direct corollary of Theorem \ref{mth1} and of the
following lemma. 

\begin{lemma} \cite[Chapter 5, \S 6.2, Theorem 1, Corollary]{Burb}\\
For any irreducible Coxeter root system $\mathcal R$ in a Euclidean
space $V$, for any $u,v \in V$
$$
\sum_{\a \in \mathcal R} \frac{(\a,u)(\a,v)}{(\a,\a)}=h (u,v),
$$
where $h$ is the Coxeter number of $\mathcal R$.
\end{lemma}

Submodules appearing in Theorem \ref{mth2} correspond to the values $c=1/m$ where $m$ is the Coxeter number of a parabolic subgroup of $W$. It follows from the description of the singular multiplicities that is the multiplicities when the polynomial representation is reducible \cite{DJO} that the multiplicity $c=k/m$ is then singular too if $k \in \N$ is coprime to $m$ (see also \cite{R}). Note however that not all the singular multiplicities have the latter form in general. For example, $c=1/9$ is singular for ${\cal E}_6$ but any parabolic subgroup of ${\cal E}_6$ has the Coxeter number at most 8. More generally the singular values $c=1/d$ where $d\in \N$ is not the Coxeter number of any parabolic subgroup of $W$ correspond to the cuspidal numbers $d$ of $W$ (see \cite{BE}). In this case any quotient of the polynomial representation of $H_c$ over its non-trivial submodule is finite dimensional \cite{BE}.

For the case of different multiplicities we define the generalized
Coxeter number $h^c=h^c_{\mathcal R}$  for the irreducible Coxeter
root system $\mathcal R$ as the coefficient of proportionality
between the following two $W$-invariant inner products
$$
\sum_{\a \in \mathcal R} \frac{c_\a(\a,u)(\a,v)}{(\a,\a)}=h^c (u,v).
$$
In the case $c=1$ we have $h^1=h$ is the usual Coxeter number. Then
Theorem~\ref{mth1} has the following reformulation.

\begin{theorem}\label{mth1'}
Let $\Gamma_0=\coprod_{i=1}^k \Gamma_i$ be the decomposition of the
subgraph $\Gamma_0 \subset \Gamma$ into the connected components.
Then the parabolic ideal $I_{\Gamma_0}$ is invariant under the
algebra $H_c$ if and only if the generalized Coxeter numbers $h^c_i$
for the Coxeter root systems determined by subgraphs $\Gamma_i$
satisfy $h^c_i=1$ for all $i=1,\ldots,k$.
\end{theorem}

\begin{remark}
It would be interesting to see if Theorem \ref{mth1'} can be established using induction and restriction functors from \cite{BE}. Let ${\cal R}_0$ be the root system with the Coxeter graph $\Gamma_0$. At the values of $c$ under consideration the corresponding rational Cherednik algebra $H_c^{{\cal R}_0}$ has trivial one-dimensional module $L$. The induced module $Ind(L)$ for $H_c^{\cal R}$ is not generally contained in the polynomial module  $\C[x]$. However we note that the modules $Ind(L)$ and $\C[x]/{I_{\Gamma_0}}$ have same support $D_{\Gamma_0}$.
\end{remark}

Our considerations allow to determine all radical ideals $I$ which
are submodules for the algebra $H_c$ in the polynomial
representation. Indeed, it is shown in \cite{BE}, \cite{Gi} that any
radical ideal $I$ must consist of polynomials vanishing on the union
of some parabolic strata $D_{\Gamma_s}$, $s=1,\ldots,L$: \beq{radI}
I=\{p\in \C[x]\mid p|_{\bigcup_{s=1}^L D_{\Gamma_s}}=0 \}. \eeq We
assume that all the strata included in the union of parabolic strata
$\bigcup D_{\Gamma_s}$ are essential in the sense that
$$
\bigcup_{s\ne j}  D_{\Gamma_s} \subsetneq \bigcup_{s=1}^L
D_{\Gamma_s}
$$
for any $1\le j \le L$.

Let $\Gamma_s=\Gamma_s^1 \bigsqcup \ldots \bigsqcup \Gamma_s^{k_s}$
be decomposition of the corresponding Coxeter graphs into the
connected components.

\begin{theorem}\label{radidth}
The radical ideal \mref{radI} is $H_c$-invariant if and only if
$h^c(\Gamma_s^i)=1$ $\forall s=1,\ldots,L; i=1,\ldots,k_s$.
\end{theorem}

The proof is similar to the proof of Theorem \ref{mth1}. In order to
derive the conditions $h^c(\Gamma_s^i)=1$ we take a polynomial $p
\in I$ in the form $p=\prod_{\pi} f_{\pi}$ where $\pi$ runs over
planes forming our union of strata: $\bigcup_{s=1}^L D_{\Gamma_s}=
\bigcup \pi$ and $f_\pi$ is a generic polynomial vanishing on $\pi$.
When $\pi$ is given as
$$
\pi = \{ x\mid (\b,x)=0, \forall \b \in \Gamma_s^j\}
$$
for $1\le s \le L, 1 \le j \le k_s$ we have $f_\pi = \sum_{\b \in
\Gamma_s^j} (x,\b) f_\b$ where $f_\b$ are generic polynomials, in
particular, $f_\b|_{\pi}$ are linearly independent. We can also
assume that $f_{\tilde \pi}|_{\pi} \ne 0$ if $\tilde \pi$ is
different from $\pi$. Then like in the proof of Theorem \ref{mth1}
the calculation of $\nabla_{\xi} f|_{\pi}$ leads to the property
$h^c(\Gamma_s^i)=1$ and hence to the Theorem.

The submodules in the polynomial representation appearing in
Theorems \ref{mth1}, \ref{mth2} correspond to the radical ideals and
to the particular singular values only.
 The description of
all singular values of parameters \cite{DJO} gives, for instance, that all the values $c=(2m-1)/2$ are
singular when $m \in \Z_+$. The next proposition describes
submodules corresponding to these singular values. When $m=1$ the
proposition is a particular case of Theorem \ref{mth2} when the
subgraph $\Gamma_0$ is a Coxeter graph $A_1$, that is consists of
one vertex, and the proof is different.

\begin{prop}\label{onedimmult}
Let $\mathcal R$ be a Coxeter root system, let $c$ be invariant
multiplicity function on $\mathcal R$. Let $S_1$ be an orbit of the
corresponding Coxeter group acting on $\mathcal R$ with the
multiplicity $c_1=c(S_1)$. Let $I$ be the ideal of polynomials
having zero of order $2m-1$ on the hyperplanes $(\a,x)=0$ for any
$\a \in S_1$. Then $I$ is $H_c$-invariant if and only if
$c_1=(2m-1)/2$.
\end{prop}
{\bf Proof.} We denote by $D_{S_1}$ the parabolic stratum
$\cup_{\a\in S_1} \{ x: (\a,x)=0 \}$. An arbitrary polynomial $p(x)$
vanishing  on $D_{S_1}$  with order $2m-1$ has the form
$$
p(x)=\prod_{\a \in S_1} (\a,x)^{2m-1} f(x)
$$
for some polynomial $f$. Let $S_2=\mathcal R \setminus S_1$ be
(possibly empty) set of roots not contained in $S_1$.


Using invariance of $\prod_{\a\in S_1} (\a,x)^{2m-1}$ with respect
to reflections $s_\b$ for $\b \in S_2$ and its anti-invariance with
respect to $s_\b$ for $\b \in S_1$ we rearrange $\nabla_\xi p(x)$ as
\begin{multline}\label{expr} \sum_{\b\in S_1} \left(
\frac{(2m-1)(\b,\xi)}{(\b,x)} \prod_{\a \in S_1} (\a,x)^{2m-1} f(x)
- 2 c_1\frac{(\b,\xi)}{(\b,x)}\prod_{\a \in S_1} (\a,x)^{2m-1}
f(x)\right)
\\
- \sum_{\b \in S_2} c_\b(\b,\xi) \prod_{\a \in S_1}  (\a,x)^{2m-1}
\frac{f(x)-s_\b f(x)}{(\b,x)}
\end{multline}
modulo elements of the ideal $I$. The last sum in \mref{expr}
belongs to the ideal $I$. The first sum in \mref{expr} belongs to
the ideal if and only if $c_1=(2m-1)/2$. Proposition is proven.

\section{Restricted Calogero-Moser systems}\label{sec2}

In this section we explain how $H_c$-invariant parabolic ideals lead
to the quantum integrable systems of Calogero-Moser type. We say
that a differential operator $L$ acting in $N$-dimensional space is
{\it quantum integrable} if there exist $N$ pairwise commuting
differential operators $L_1=L, \ldots, L_N$ so that $L_i$ are
algebraically independent.

 Consider the parabolic ideal $I$
consisting of polynomials vanishing on the parabolic stratum $D$
which is the $W$-orbit of the subspace $\pi$. Assume that ideal $I$
is a submodule for the rational Cherednik algebra $H_c$. Let $p=\bar
p|_D$ be a restriction to the stratum of a polynomial $\bar p$
defined in the whole space $V$.
It follows from the invariance of the ideal $I$ that for any Dunkl
operator $\nabla_\xi$ the result of restriction $\nabla_\xi \bar
p|_D$ does not depend on the extension $\bar p$ but depends on $p$
only. Therefore the restricted Dunkl operators $\nabla_\xi|_D$ are
correctly defined.

Moreover, the analysis of the property of the ideal $I$ to be
invariant in the proof of Theorem \ref{mth1} implies that the
restricted operators $\nabla_\xi|_D$ are defined correctly in the
locally analytic settings. Namely, let point $x_0 \in\pi$ be
generic, let $U \ni x_0$ be its small neighborhood, $U \subset \pi$.
Consider its $W$-orbit $U^W=\cup_{w \in W} w(U)$. Let $f$ be a union
of analytic germs defined in the neighborhoods $U^W$, let $\bar f$
be analytic extension of these germs to $\tilde U^W=\cup_{w \in W}
w(\tilde U)$ where $\tilde U \supset U$ is a small neighborhood of
$x_0$ in the space $V$. Then the result of the restriction
$\nabla_\xi \bar f|_D$ does not depend on the locally analytic
extension $\bar f$ but depends on $f$ only.

Consider now the space $\cal L$ of $W$-invariant union of germs $f$
defined in $U^W$. So $f$ is determined by its values $f_\pi=f|_U$ in
the neighborhood $U\subset \pi$. The invariant combinations of Dunkl
operators $\sigma (\nabla)=\sigma(\nabla_1, \ldots, \nabla_N)$,
where  $\sigma(x_1, \ldots, x_N)\in \C[x_1,\ldots,x_N]^W$, act in
the space $\mathcal L$. We denote by $\sigma(\nabla)^{Res_\pi}$ the
operator which maps $f_\pi$ to the result of the restriction
$\sigma(\nabla)\bar f|_U$
 on to
the neighborhood $U \subset \pi$ of the application of the operator
$\sigma(\nabla)$ to any $W$-invariant extension $\bar f$ of $f$ from
$U^W$ to $\tilde U^W$ (c.f. \cite{Heckman}).



\begin{theorem}\label{mth3}
Assume that a stratum $D$ defines $H_c$-invariant parabolic ideal.
Then the operator $\sum_{i=1}^N \nabla_i^2$ restricted to the
$W$-invariant functions on $D$ has the generalized Calogero-Moser
form \beq{CMres} {\left(\sum_{i=1}^N \nabla_i^2\right)}^{Res_\pi}=
\Delta_y -\sum_{\nad{\a\in \mathcal R_+}{\hat\a\ne 0}} \frac{2
c_\a}{(\hat\a,y)} \p_{\hat\a}, \eeq where $y=(y_1,\ldots,y_n)$ are
orthonormal coordinates on the plane $\pi$, $ \Delta_y =
\frac{\p^2}{\p y_1^2}+\ldots + \frac{\p^2}{\p y_n^2}$, vector $\hat
\a$ is orthogonal projection of vector $\a$ onto $\pi$.

For any polynomials $\sigma, \tau \in \C[x]^W$ the restrictions
$\sigma(\nabla)^{Res_\pi}$, $\tau(\nabla)^{Res_\pi}$ are commuting
differential operators in the space $\pi$, in particular, operator
\mref{CMres} is quantum integrable.
\end{theorem}
{\bf Proof.} 
The operator
$$
H=\sum_{i=1}^N \nabla_i^2
$$ can be expanded as
\beq{CM} H = \Delta - \sum_{\a \in \mathcal R_+} \frac{2 c_\a
}{(\a,x)}\p_\a + \sum_{\a \in \mathcal R_+} \frac{c_\a
(\a,\a)(1-s_{\a})}{(\a,x)^2}. \eeq

Consider $f$ which is a $W$-invariant analytic function defined in
the neighborhoods $w(U) \subset w(\pi)$ of the $W$-orbit of the
generic point $x_0 \in\pi$. Consider now invariant analytic
extension $\bar f$ of the function $f$ to the union of neighborhoods
$w(\tilde U)$ where $\tilde U \supset U$ is a neighborhood of $x_0$
in $\C^N$, $w \in W$. We are going to apply the operator $H$ to the
function $\bar f$.
 The assumption of
the theorem says that the result of the restriction of $H \bar f$
onto $\pi$ does not depend on the extension $\bar f$ and it is
determined by $f_\pi=\bar f|_U$ only. So we may choose $\bar f$ to
be constant along the normal directions to $\pi$. Then
$$
\p_\a \bar f = \p_{\hat\a} \bar f, \quad \Delta \bar f = \Delta_y
\bar f,
$$
and also $(\a,x) = (\hat\a,x)$ when $x\in \pi$. Since function $\bar
f$ is $W$-invariant, the last sum in \mref{CM} disappears and the
Calogero-Moser operator \mref{CM} takes the form \mref{CMres}. This
proves the first part of the theorem.

The second statement follows from the commutativity of the Dunkl
operators $[\nabla_\xi, \nabla_\eta]=0$ and from the fact that the
operators $\sigma(\nabla)$ preserve the space $\mathcal L$ of
$W$-invariant germs. The highest terms of the restrictions $\sigma(\nabla)^{Res_\pi}$ are obtained by the restriction of the highest term of the operator $\sigma(\nabla)^{Res}$ onto the plane $\pi$. Since the stratum $D$ considered in the space of orbits $\C^N/W$ is a (singular) variety of dimension $n$ there exist polynomials $\sigma_1=x_1^2+\ldots+x_N^2,\sigma_2,\ldots,\sigma_n \in \C[x]^W$ such that their restrictions on $D$ are algebraically independent. The corresponding differential operators are algebraically independent as well, hence the operator \mref{CM} is quantum integrable. This completes the proof of the theorem.

We  note that the CM system itself can also be restricted to the
stratum $D$ considered inside the orbit space if the corresponding
parabolic ideal is $H_c$-invariant. Also it is known that in the
orbit space the CM system becomes algebraic \cite{Turb} so we have
restriction of the non-singular differential operator to a
subvariety.

Specific choice of invariants in Theorem \ref{mth3} leading to a collection of algebraically independent differential operators on $\pi$ depends of course on the particular Coxeter group $W$ and the stratum. In the case of classical Coxeter groups one can always take the Newton sums as such invariant polynomials. More exactly let $\sigma_k = \sum_{i=1}^N x_i^k$. Let $\pi$ be an intersection of mirrors of the group $W$ of dimension $n$. In the case $W={\mathcal A}_{N-1}$ the polynomials
$\sigma_1|_\pi,\ldots,\sigma_n|_\pi$ and the corresponding differential operators on $\pi$ are algebraically independent. This follows from the explicit form of the $H_c$-invariant strata and the fact that deformed Newton sums $\hat \sigma_k = \sum_{i=1}^{n_1} y_i^k + \kappa \sum_{i=1}^{n_2} z_i^k$ are algebraically independent for $k=1,\ldots, n=n_1+n_2$ if $\kappa\in \N$ (this in turn follows from  \cite[Proposition 4]{SV1}).
In the case $W={\mathcal B}_N$ or $W={\mathcal D}_N$ the polynomials $\sigma_{2i}|_\pi$, $i=1,\ldots,n$  and  the corresponding differential operators are algebraically independent by similar reasons.

So far we were using  the ``radial normalization" \mref{CMres} of
the generalized CM systems. The restricted operators are also gauge
equivalent to the operators in the ``potential normalization". More
exactly, we have the following property of the generalized CM
systems related to arbitrary  parabolic stratum.

\begin{prop}
Let $\pi\subset V$ be an intersection of mirrors \beq{pidef2} \pi =
\{x\in V \mid (\b,x)=0 \quad \forall \b\in \Gamma_0^v \}, \eeq
corresponding to a Coxeter subgraph $\Gamma_0 \subset \Gamma$. Let
$\hat u$ denote the orthogonal projection of a vector $u \in V$ onto
the space $\pi$.
Consider the decomposition $\hat {\mathcal R}_+ = R_1 \sqcup \ldots \sqcup R_k$ such that for any two vectors $u,v \in R_i$ we have collinearity
$u \sim v$, and for any two non-zero vectors $u \in R_i$, $v \in R_j$ one has $u \nsim v$ when $i \ne j$.
 Choose a nonzero element $\hat \gamma_i \in R_i$ for any $i=1,\ldots, k$. Define $c_i = \sum_{\nad{\a \in {\mathcal R}_+}{\hat \a \in R_i \setminus 0}} c_\a$.

Then \beq{gauge} f^{-1} (\Delta-\sum_{\nad{\a\in
{\mathcal R}_+}{\hat\a\ne 0}} \frac{2 c_\a}{(\hat\a,x)}\p_{\hat\a})
f = \Delta - \sum_{i=1}^k
\frac{c_i (c_i+1)(\hat\gamma_i,\hat\gamma_i)}{(\hat\gamma_i,x)^2} \eeq  where $f=
\prod_{\nad{\a \in {\mathcal R}_+}{\hat\a\ne 0}} (\hat
\a,x)^{c_\a}$, and $\Delta$ is Laplacian on $\pi$.
\end{prop}
{\bf Proof.} The gauge property \mref{gauge} is equivalent to the
following series of identities for all $\a \in \mathcal R$:
\beq{togd} \sum_{\nad{\b \in {\mathcal R}}{\hat\b \nsim \hat\a}}
\frac{c_\b (\hat \a, \hat \b)}{(\hat \b,x)}=0 \eeq when $x \in \pi$
and $(\hat \a,x)=0$. To establish identity \mref{togd} consider the
set $S\subset \mathcal R$ of the roots $\b$ such that $\hat \b$ is
not proportional to $\hat \a$, we also assume that $\hat \a \ne 0$.
Consider the action on $S$ of the group $W_0$ which is generated by
reflections at the simple roots $\Gamma_0^v$ and by the reflection
$s_\a$. We have decomposition of the set $S$ as the union of
$W_0$-orbits $S= \mathcal O_1 \cup \ldots \cup \mathcal O_k$. For
each $i=1,\ldots,k$ we claim that \beq{togd2} \sum_{\b \in \mathcal
O_i} \frac{(\hat \a, \hat \b)}{(\hat \b,x)}=0 \eeq when $x \in \pi$
and $(\hat \a,x)=0$. Indeed, for vectors $x$ under consideration, we
have $(\hat \b_1,x)=(\hat \b_2,x)$ for any roots $\b_1,\b_2 \in
\mathcal O_i$. Also vector $b:=\sum_{\b \in \mathcal O_i} \b$
satisfies \beq{togd3} ( \alpha, b)= (\gamma, b)=0 \eeq for any
$\gamma \in \Gamma_0^v$ because of invariance $s_\a b = b = s_\g b$.
Property \mref{togd3} implies $\hat b = b$ and then the identities
\mref{togd2} hold which in turn imply \mref{togd}. The Proposition
is proven.

\section{Parabolic ideals and CM systems from classical Coxeter
groups}\label{sec3}

We are going now to apply the described method of restriction to
derive particular integrable systems of Calogero-Moser type. In this
section we deal with the rational Cherednik algebras for the
classical root systems, exceptional Coxeter root systems are dealt
with in Section \ref{sec5}.

\subsection{${\mathcal A}$-series}
Consider the root system $\mathcal A_{N-1}\subset \C^N$ given by the
collection of vectors $e_i-e_j$, $1 \le i<j \le N$. Theorems
\ref{mth1}-\ref{mth1'} remain true when the rank of the root system
is less than the dimension of the ambient space $V$.
 As it follows from
Theorem \ref{mth2} the parabolic strata which define invariant
ideals for the rational Cherednik algebra for the ${\mathcal
A}_{N-1}$ root system must have the form \beq{diskrmk}
D_{m,k}=\bigcup_{w\in S_N} w(\pi_{m,k}) \eeq where the plane
$\pi_{m,k}$ is given by the equations
\begin{multline}\label{pimk}
x_1=x_2=\ldots=x_{k},\quad x_{k+1}=x_{k+2}=\ldots=x_{2 k},\\
\ldots\\
x_{(m-1)k+1}= x_{(m-1)k+2}= \ldots =x_{mk}.
\end{multline}
Here $m \ge 1, k>1$ are integer such that $mk \le N$, the
corresponding parameter $c=1/k$.

\begin{prop}\label{prop3}
For the root system ${\mathcal A}_{N-1}$, the parabolic stratum
\mref{diskrmk}, \mref{pimk} and the multiplicity $c=1/k$ the
restricted Calogero-Moser operator takes the form
\begin{multline}\label{amnop}
H=\left(\sum_{i=1}^N \nabla_i^2\right)^{Res_{\pi_{m,k}}} =
 \Delta - 2 k \sum_{1\le i < j\le m}
\frac{\p_i-\p_j}{y_i-y_ j} -
\\
 \frac{2}{k} \sum_{m+1 \le
i < j \le m+ n}\frac{\p_i-\p_j}{y_i-y_ j} - 2 \sum_{\nad{1\le i \le
m}{m+1 \le j\le m+ n}} \frac{\p_i-\sqrt{k}\p_j}{y_i-\sqrt{k} y_ j},
\end{multline}
where $\Delta=\sum_{i=1}^{m+n}\p_i^2$, $\p_i=\frac{\p}{\p y_i}$, and
$n=N-mk$.
\end{prop}
{\bf Proof.} One way to deduce this statement is to consider the
projection of the root system ${\mathcal A}_{N-1}$ to the plane
$\pi_{m,k}$ as Theorem \ref{mth3} directs. Alternatively, we may
introduce the following change of coordinates. For every block of
colliding coordinates $x_{jk+1}=\ldots=x_{jk+k}$ we define new
coordinates $z_{jk+1},\ldots,z_{jk+k}$ such that
$$
z_{jk+1}=\frac{x_{jk+1}+\ldots+x_{jk+k}}{\sqrt{k}},
z_{jk+2}=\frac{x_{jk+1}-x_{jk+2}}{\sqrt{k}},
z_{jk+3}=\frac{x_{jk+2}-x_{jk+3}}{\sqrt{k}}, \ldots,
$$
$$
z_{jk+k}=\frac{x_{jk+k-1}-x_{jk+k}}{\sqrt{k}}
$$
where $j=0,\ldots,m-1$. The remaining coordinates are not changed:
$z_i=x_i$ for $i>mk$. Then on the plane $\pi_{m,k}$ we have
$x_{jk+s}=\frac{1}{\sqrt{k}}z_{jk+1}$, $s=1,\ldots,k$. Also
$\frac{\p}{\p x_{jk+s}}=\frac{1}{\sqrt{k}} \frac{\p}{\p z_{jk+1}}$
when acting on functions $\bar f$ which are constant along the
directions orthogonal to the plane $\pi_{m,k}$. One gets operator
\mref{amnop} after renaming surviving $z$-coordinates into
$y$-coordinates.

The operator \mref{amnop} in the case $n=1$ appeared first in the
work \cite{CFV1}  where its (algebraic) integrability was
established. For arbitrary $m,n$ and $k$ the integrability of the
trigonometric version of the operator \mref{amnop} was established
in \cite{SV1} using explicit calculations. In the work \cite{SV2}
integrability for arbitrary $m,n,k$ was established by obtaining the
operator as a restriction of the trigonometric Calogero-Moser
operator acting in the space of symmetric functions of infinitely
many variables to the generalized discriminants.

\subsection{$\mathcal B$-series}

Consider the group ${\mathcal B}_N$ generated by reflections at
$x_i=\pm x_j$ for $1\le i<j \le N$, and at $x_i=0$ for
$i=1,\ldots,N$. Let $c_2$ be the multiplicity of the roots $e_i$,
and let $c_1$ be the multiplicity of the roots $e_i \pm e_j$. There
are parabolic strata corresponding to the Coxeter subgraphs of the
type ${\mathcal A}_{k-1}\times\ldots \times {\mathcal A}_{k-1}$,
$k>1$, where the number of subsystems ${\mathcal A}_{k-1}$ is $m$,
$mk\le N$. The corresponding plane $\pi_{m,k}\subset \C^N$ is given
by the equations
\begin{multline}\label{pimkb}
x_1=x_2=\ldots=x_{k},\quad x_{k+1}= x_{k+2}=
\ldots= x_{2 k},\\
\ldots\\
x_{(m-1)k+1}= x_{(m-1)k+2}= \ldots = x_{mk}.
\end{multline}
The corresponding parabolic stratum $D_{m,k}\subset \C^N$  is the
orbit of the plane $\pi_{m,k}$ under the group ${\mathcal B}_N$:
\beq{stratumBN1} D_{m,k}=\bigcup_{w\in {\mathcal B}_N} w(\pi_{m,k}).
\eeq Theorem \ref{mth2} and easy calculations imply the following
\begin{prop}\label{BN1}
For the root system ${\mathcal B}_{N}$ with the multiplicity
$c_1=1/k$
 the
parabolic ideal
$$I_{m,k}=\{f\in \C[x_1,\ldots,x_N] \ | \ f|_{D_{m,k}}=0\}$$
corresponding to the stratum \mref{pimkb}, \mref{stratumBN1} is
invariant under the rational Cherednik algebra $H_c^{{\mathcal
B}_N}$. The restricted Calogero-Moser operator takes the form

\begin{multline}\label{bn1} H=\Delta-\sum_{\nad{i,j=1}{i\ne j}}^m
\left(\frac{2k(\p_i-\p_j)}{y_i-y_j}+\frac{2k(\p_i+\p_j)}{y_i+y_j}
\right)-\sum_{i=1}^m\sum_{j=m+1}^{m+n}
\left(\frac{2(\p_i-\sqrt{k}\p_j)}{y_i-\sqrt{k}y_j}+\frac{2(\p_i+\sqrt{k}\p_j)}{y_i+\sqrt{k}y_j}
\right)
\\ - \sum_{\nad{i,j=m+1}{i\ne j}}^{m+n}
\left(\frac{2k^{-1}(\p_i-\p_j)}{y_i-y_j}+\frac{2k^{-1}(\p_i+\p_j)}{y_i+y_j}
\right)-\sum_{i=1}^m \frac{(kq+k-1)\p_i}{y_i}-\sum_{i=m+1}^{m+n}
\frac{q\p_i}{y_i},
\end{multline}
where $\Delta=\sum_{i=1}^{m+n} \p_i^2$, $q=2c_2 \in \C$, and
$n=N-mk$.
\end{prop}

For $n=1$ the first sum in the second line in \mref{bn1} disappears
and the integrability of this operator was established in
\cite{CFV3} where it corresponded to the configuration
${\mathcal C}_{m+1}(\frac{kq+k-1}2,\frac{q}2)$. In the case $n>1$ the
integrability of the trigonometric version of \mref{bn1} was shown
in \cite{SV1},\cite{SV3}.

Consider the plane $\pi_{l}\subset \C^N$ given by the equations
$$ x_1=x_2=\ldots=x_{l}=0. $$ Consider the corresponding
parabolic stratum $D_{l}$ which is the orbit of the plane $\pi_{l}$
under the group ${\mathcal B}_N$: \beq{stratumBN2} D_{l}=\bigcup_{i_1<\ldots
<i_l} \{x\mid x_{i_1}=\ldots=x_{i_l}=0\}. \eeq Define the ideal
$I_{l}=\{f\in \C[x_1,\ldots,x_N]  \mid f|_{D_{l}}=0\}$.

\begin{prop}\label{BN2}
The Dunkl operators $\nabla_i$ preserve the ideal $I_l$ if and only
if $2(l-1)c_1+2c_2=1$. The corresponding restricted Calogero-Moser
operator is the CM operator constructed by the root system
${\mathcal B}_{N-l}$.
\end{prop}
{\bf Proof.} It is easy to check that the generalized Coxeter number
$h_c$ for the root system ${\mathcal B}_l$ with the multiplicities
$c(e_i \pm e_j)=c_1$, $c(e_i)=c_2$ equals $h_c=2c_1(l-1)+2c_2$ . The
statement follows from Theorem \ref{mth1'}.

\begin{theorem}
All possible parabolic strata defining
$H_c^{\mathcal{B}_N}$-invariant ideals are either the strata
$D_{m,k}$ defined by \mref{stratumBN1} (with $c_1=1/k$) or strata
$D_l$ defined by \mref{stratumBN2} (with $2(l-1)c_1+2c_2=1$) or
their intersection.
In the latter case the stratum is the ${\mathcal B}_N$-orbit of the subspace
\begin{multline}
x_1=x_2=\ldots=x_{k},\quad x_{k+1}= x_{k+2}=
\ldots= x_{2 k},\\
\ldots\\
x_{(m-1)k+1}= x_{(m-1)k+2}= \ldots = x_{mk},
\\ x_{mk+1}=\ldots=x_{mk+l}=0,
\end{multline}
with $c_1=1/k, c_2=\frac12 -\frac{l-1}{k}$, and $k>1$, $l \ge 1$,
$mk+l \le N$. In this case the corresponding restricted
Calogero-Moser operator is operator \mref{bn1} where
$q=\frac{k+2l+2}{k}$.
\end{theorem}
The proof follows from the general structure of the parabolic strata
for the ${\mathcal B}_N$ group, Theorem \ref{mth1'} and Propositions
\ref{BN1}, \ref{BN2}.

\subsection{$\mathcal D$-series}

Consider the group ${\mathcal D}_N$ generated by reflections at
$x_i=\pm x_j$ for $1\le i<j \le N$.
 Take two integer parameters $m>0,k>1$ such that
$m k \le N$. Consider the plane $\pi_{m,k}^\varepsilon\subset \C^N$
given by the equations
\begin{multline}\label{pimkd}
x_1=x_2=\ldots=x_{k},\quad x_{k+1}= x_{k+2}=
\ldots= x_{2 k},\\
\ldots\\
x_{(m-1)k+1}= x_{(m-1)k+2}= \ldots = \varepsilon x_{mk}.
\end{multline}
Here $\varepsilon=1$ except the case when $k$ is even and $N=mk$. In
the latter case $\varepsilon=\pm 1$.
 Consider the corresponding parabolic stratum
$D_{m,k}^\varepsilon\subset \C^N$ which is the orbit of the plane
$\pi_{m,k}^\varepsilon$ under the group ${\mathcal D}_N$:
\beq{stratumDN1} D_{m,k}^\varepsilon=\bigcup_{w\in {\mathcal D}_N}
w(\pi_{m,k}^\varepsilon). \eeq This describes all possible strata in
${\mathcal D}_N$ with the Coxeter graphs ${\mathcal
A}_{k-1}\times\ldots \times {\mathcal A}_{k-1}$ (see \cite{OS}).

Theorem \ref{mth2} and easy calculations imply the following
\begin{prop}
For the root system ${\mathcal D}_{N}$, the parabolic stratum
\mref{pimkd}, \mref{stratumDN1} and the multiplicity $c=1/k$ the
corresponding ideal $$I_{m,k}^\ve=\{f\in \C[x_1,\ldots,x_N] \ | \
f|_{D_{m,k}^\varepsilon}=0\}$$ is invariant under the corresponding
rational Cherednik algebra $H_c$. The restricted Calogero-Moser
operator has the form \mref{bn1} with $q=0$.
\end{prop}

Another type of the parabolic strata for the ${\mathcal D}_N$ group
has the Coxeter graph of type ${\mathcal D}_p$, $1<p<N$. Consider
the plane $\pi_{p}\subset \C^N$ given by the equations
$$ x_1=x_2=\ldots=x_{p}=0. $$ The corresponding
parabolic stratum $ D_{p}$ is the orbit of the plane $\pi_{p}$ under
the group ${\mathcal D}_N$: \beq{stratumDN2}
D_{p}=\bigcup_{i_1<\ldots <i_p} \{x\mid x_{i_1}=\ldots=x_{i_p}=0\}.
\eeq Theorem \ref{mth2} implies that the ideal $I_{p}=\{f\in
\C[x_1,\ldots,x_N] \mid f|_{D_{p}}=0\}$ is invariant under the rational
Cherednik algebra if and only if $c=\frac{1}{2(p-1)}$. The
corresponding restricted Calogero-Moser operator is the CM operator
constructed by the root system ${\mathcal B}_{N-p}$.

The next theorem is a corollary of  Theorem \ref{mth2}, previous
considerations and calculations of the restricted root systems.

\begin{theorem}
For the root system ${\mathcal D}_N$ all possible parabolic strata
defining $H_c^{\mathcal D_N}$-invariant ideals  are either the
strata $D_{m,k}^\ve$ defined by \mref{stratumDN1} (with $c=1/k$) or
strata $D_p$ defined by \mref{stratumDN2} (with
$c=\frac{1}{2(p-1)}$) or their intersection.
In the latter case the stratum is the ${\mathcal D}_N$-orbit of the
plane
\begin{multline}\label{dnideal}
x_1=x_2=\ldots=x_{k},\quad x_{k+1}= x_{k+2}=
\ldots=  x_{2 k},\\
\ldots\\
x_{(m-1)k+1}= x_{(m-1)k+2}=
\ldots =  x_{mk}, \\
x_{mk+1}=\ldots=x_{mk+\frac{k}{2}+1}=0,
\end{multline}
where $k\ge 2$ is even, and $c=1/k$. In this case the restricted
Calogero-Moser operator is the operator \mref{bn1} with
$q=\frac{2(k+2)}{k}$.\end{theorem}

\section{Generalized CM systems with quadratic
potential}\label{sec4}

In this section we show that generalized CM operators \mref{amnop},
\mref{bn1} obtained from the parabolic strata for the classical
Coxeter groups remain integrable when the term $\omega^2
\sum_{i=1}^{m+n} y_i^2$ is added to the Hamiltonians $H$.

First we recall integrability of the Calogero-Moser systems with
square potential following \cite{Poly}. Let $\nabla_i$ be Dunkl
operator in the basis direction $e_i$ for the symmetric group $S_N$,
$$
\nabla_i=\p_i-c \sum_{\nad{j=1}{j\ne i}}^N
\frac1{x_i-x_j}(1-s_{ij}).
$$
  Define now the
operators $\nabla_i^{+}, \nabla_i^{-}$ as
$$
\nabla_i^{\pm} = \nabla_i \pm \omega x_i,
$$
where $\omega \in \C$ is a parameter. A combination of operators
$\nabla_i^\pm$ gives CM system with quadratic potential. Namely,
because of commutation relations \beq{xnablasame} [x_i,
\nabla_i]=-1+ c \sum_{\nad{j=1}{j\ne i}}^N s_{ij},  \eeq we have
$$
\sum_{i=1}^N \nabla_i^+ \nabla_i^-= \sum_{i=1}^N \nabla_i^2
-\omega^2 \sum_{i=1}^N x_i^2 -\omega N + 2\omega c \sum_{{i<j}}^N
s_{ij}.
$$
The last expression becomes CM Hamiltonian with the added quadratic
terms and constant when restricted to the space of symmetric
functions. Denoting the operation of restriction to invariants by
$Res$ (c.f. \cite{Heckman}) we recall the following result.

\begin{theorem}\cite{Poly}\label{polyth}
The CM Hamiltonian with quadratic potential
$$
H=\Delta-\sum_{i<j}^N \frac{2 c (\p_i-\p_j)}{(x_i-x_j)} -\omega^2
\sum_{i=1}^N x_i^2
$$
can be obtained as
$$
H=\left(\sum_{i=1}^N \nabla_i^+ \nabla_i^-\right)^{Res}+ \omega N-
\omega c N(N-1).
$$
The differential operators $\left(\sum_{i=1}^N \left(\nabla_i^+
\nabla_i^-\right)^k \right)^{Res}$ for $k \in \Z_+$ pairwise
commute, so $H$ is quantum integrable.
\end{theorem}
{\bf Proof.} Commutativity of the differential operators follows
from the fact that the combinations $\sum_{i=1}^N (\nabla_i^+
\nabla_i^-)^k$ preserve the space of symmetric polynomials and from
the following commutativity of the combinations of Dunkl operators
\beq{hcom} [\sum_{i=1}^N h_i^k, \sum_{i=1}^N h_i^l]=0, \eeq where
$h_i=\nabla_i^+ \nabla_i^-$ and $k,l \in \Z_+$. To establish
\mref{hcom} we note first that from the relations $[x_i,
\nabla_j]=-c s_{ij}$ valid for $i\ne j$ it follows that
$[h_i,h_j]=2\omega c (h_i-h_j) s_{ij}$. Then by induction in $k$ it
is easy to deduce that
$$
[h_i^k, h_j] = 2\omega c (h_i^k-h_j^k)s_{ij},
$$
and then by induction in $l$ one obtains that
$$
[h_i^k, h_j^l] = 2\omega c \sum_{t=1}^l
(h_j^{t-1}h_i^{k+l-t}-h_j^{k+l-t}h_i^{t-1})s_{ij}.
$$
The last formula implies by induction in $l$ that $[h_i^k, h_j^l]$
is antisymmetric with respect to permutation $h_i$ and $h_j$.
Therefore $[h_i^k, h_j^l]+[h_j^k, h_i^l]=0$, hence commutativity of
the operators \mref{hcom} holds, and the theorem follows.

Similar arguments allow to establish quantum integrability of the
Calogero-Moser systems with quadratic potential in the case of other
classical Coxeter root systems \cite{Chalykhq}.

Indeed, let now $ \nabla_i^{\pm} = \nabla_i \pm \omega x_i, $ where
$\nabla_i$ is the Dunkl operator in the direction $e_i$ for the root
system ${\mathcal B}_N$ or ${\mathcal D}_N$, $i=1,\ldots,N$. The
operators $\sum_{i=1}^N h_i^k$ where $h_i=\nabla_i^+ \nabla_i^-$
preserve the spaces of the corresponding $\mathcal B_N$ or $\mathcal
D_N$ invariants. Commutativity of these combinations for different
$k$ can be established similar to the proof of commutativity
\mref{hcom} in the proof of Theorem \ref{polyth}. Indeed, the
commutation relations need to be modified as follows
$$[x_i, \nabla_j]=-c (s_{ij}-s_{ij}^+), \, [h_i,h_j]=2\omega c (h_i-h_j)
(s_{ij}+s_{ij}^+), \, [h_i^k, h_j] = 2\omega c
(h_i^k-h_j^k)(s_{ij}+s_{ij}^+),
$$
where $1\le i <j \le N$, and then
$$
[h_i^k, h_j^l] = 2\omega c \sum_{t=1}^l
(h_j^{t-1}h_i^{k+l-t}-h_j^{k+l-t}h_i^{t-1})(s_{ij}+s_{ij}^+),
$$
where $s_{ij}^+$ is the reflection at the hyperplane $x_i+x_j=0$.

Finally, the differential operators $(\sum_{i=1}^N \nabla_i^+
\nabla_i^-)^{Res}$ coincides with the corresponding ${\mathcal B}_N$
or ${\mathcal D}_N$ CM operator with quadratic potential up to a
constant, if $Res$ denotes the restriction to ${\mathcal B}_N$ or
${\mathcal D}_N$-invariants respectively (c.f. \cite{Heckman}).

 Now we are ready to obtain
generalized quantum Calogero-Moser systems with the additional
quadratic potential.

\begin{theorem}\label{th8}
The operators  $H-\omega^2 \sum_{i=1}^{m+n} y_i^2$ are quantum
integrable if $H$ is given either by the formula \mref{amnop} or by
the formula \mref{bn1}.
\end{theorem}
In this theorem we assume that $k \in \Z$ as it also happens in
Propositions \ref{prop3}, \ref{BN1}, although we expect Theorem
\ref{th8} to be true for any $k$. When $H$ is given by \mref{amnop}
with $n=1$ Theorem \ref{th8} is already established in \cite{CFV2}
for any $k\in \C$. We also note that some eigenfunctions of the
operators from the theorem were already investigated in \cite{HL,
H3}.

The proof of the theorem follows from the fact that if a parabolic
ideal $I$ corresponding to the Coxeter orbit of a linear subspace
$\pi$ is preserved under the action of the Dunkl operators
$\nabla_i$, then $I$ is also preserved under the action of the
operators $\nabla_i^\pm$. Indeed, consider the combinations
$$
L_k=\sum_{i=1}^N h_i^k
$$
where $h_i=\nabla_i^+ \nabla_i^-$ is a product of Dunkl operators
for the Coxeter groups $\mathcal A_{N-1}$ or $\mathcal B_N$. The
operators $L_k$ are invariant under the action of the corresponding
classical Coxeter group. The restrictions $L_k^{res_\pi}$ are
commuting differential operators on the subspace $\pi$. For the
Coxeter group $\mathcal A_{N-1}$ with multiplicity $c=k^{-1}$ and
the subspace $\pi$ given by the equations \mref{pimk} the operator
$$
L_1^{Res_\pi}=H-\omega^2\sum_{i=1}^{m+n} y_i^2 -\omega N + \omega
k^{-1} N(N-1),
$$
where $H$ is given by \mref{amnop}. For the Coxeter group $\mathcal
B_{N}$ with multiplicities $c(e_i \pm e_j)=k^{-1}$, $c(e_i)=q/2$ and
the subspace $\pi$ given by the equations \mref{pimkb} the operator
$$
L_1^{Res_\pi}=H-\omega^2\sum_{i=1}^{m+n} y_i^2 -\omega N + 2\omega
k^{-1} N(N-1) +\omega q N,
$$
where $H$ is given by \mref{bn1}.

\section{Generalized CM systems from exceptional Coxeter groups}\label{sec5}

First we consider invariant ideals for the Coxeter root system
$\mathcal F_4$. We use Theorem \ref{mth1'} and the fact that the
parabolic strata are given by non-isomorphic Coxeter subgraphs and
by isomorphic subgraphs if they are different as Dynkin subgraphs of
the root system $\mathcal F_4$ (\cite{OS}). We also use that for the
root system ${\mathcal A}_n$ with multiplicity $c$ we have
$h^c=(n+1)c$, and for the root system ${\mathcal B}_n$ with
multiplicities $c(e_i \pm e_j)=c_1$, $c(e_i)=c_2$ we have
$h^c=2(n-1)c_1+2c_2$.

The action of the group on the root system $\mathcal F_4$ has two
orbits with multiplicities $c_1$, $c_2$. There are two strata
corresponding to the subgraphs of type $\mathcal A_1$, the
corresponding ideals are invariant when the corresponding
multiplicity $c_1$ or $c_2$ equals $1/2$. There are two strata
corresponding to subgraphs of type $\mathcal A_2$. They are
invariant iff the corresponding multiplicity $c|_{{\mathcal
A}_2}=1/3$. There is one stratum of the type $\mathcal A_1 \times
\mathcal A_1$, it is invariant iff $c_1=c_2=1/2$. The stratum of
type $\mathcal B_2$ is invariant iff $c_1+c_2=1/2$. There are two
strata of type $\mathcal B_3$, they are invariant iff the
corresponding generalized Coxeter number $h^c_{{\mathcal B_3}}=1$
which gives $4c_1+2 c_2=1$ or $4c_2+2 c_1=1$ respectively. There are
two strata of type $\mathcal A_1 \times \mathcal A_2$, they are
invariant iff $c|_{{\mathcal A}_2}=1/3$ and $c|_{{\mathcal
A}_1}=1/2$.

The restricted Calogero-Moser systems correspond to the root system
$\mathcal G_2$ for the strata $\mathcal A_2$, and the restricted
system corresponds to the root system $\mathcal B_2$ for the stratum
$\mathcal B_2$.

The restricted Calogero-Moser systems corresponding to the strata
$\mathcal A_1$ give equivalent new non-Coxeter one-parametric
families  of integrable systems in dimension 3.

\begin{theorem}\label{theoremf4a1}
The restricted $\mathcal F_4$ CM system for the stratum $\mathcal
A_1$ has the form \beq{f4a1} H=\Delta- \sum_{i=1}^3
\frac{(4c+1)\p_i}{x_i}-\sum_{\nad{i,j=1}{{\nad{i<j}{ }}}}^3
\frac{2c(\p_i\pm \p_j)}{x_i \pm x_j}- \sum \frac{2(\p_1 \pm \p_2 \pm
\p_3)}{x_1\pm x_2 \pm x_3}, \eeq where $c \in \C$ and summations run
over arbitrary choices of signs. In particular, operator \mref{f4a1}
is quantum integrable.
\end{theorem}

In the case of the $\mathcal A_1 \times \mathcal A_1$ strata the
restricted integrable CM system has the Hamiltonian
\begin{multline}\label{f4a1a1} H=\p^2_x+\p^2_y - \frac{2m}{x}\p_x-\frac{2m}{y}\p_y
- \frac{2n}{x+y}(\p_x+\p_y) - \frac{2n}{x-y}(\p_x-\p_y)- \\
\frac{2(\p_x+\a \p_y)}{x+\a y}-\frac{2(\p_x-\a \p_y)}{x-\a y}-
\frac{2(\p_x+\a^{-1} \p_y)}{x+\a^{-1} y}-\frac{2(\p_x-\a^{-1}
\p_y)}{x-\a^{-1} y}
\end{multline} where $m=7/2$, $n=0$ and $\a=\sqrt{2}$.

\begin{prop}\label{propos} {\it The system \mref{f4a1a1} is quantum integrable for any
$m$, $n$ and $\a=\pm(\sqrt{2n+1} \pm \sqrt{2(m+n+1)})/\sqrt{2m+1}$}.
\end{prop}
Proof follows from the fact that operator \mref{f4a1a1} in the
potential form satisfies the locus conditions (\cite{CFV3}) when
$m,n\in \Z$. There is a commuting operator to \mref{f4a1a1} with the
highest symbol $(\a^2\xi_1^2-\xi_2^2)^3 (\xi_1^2-\a^2\xi_2^2)^3$.




Now we give a complete list of the sets $A$ of vectors $\a$ with
multiplicities $m_{\a}$ for the  generalized CM systems,
which are obtained as restrictions of the exceptional Coxeter root
systems of $\mathcal E$ and $\mathcal H$ type. The corresponding
operators
$$
\Delta - \sum_{\a \in  A} \frac{2 m_\a}{(\a,x)}\p_\a,\quad \Delta -
\sum_{\a \in  A} \frac{m_\a(m_\a+1)(\a,\a)}{(\a,x)^2}
$$
are quantum integrable. These systems are labeled by a pair
$(\Gamma,\Gamma_0)$, where $\Gamma$ is a Coxeter graph of $\mathcal
E$ or $\mathcal H$ type, and $\Gamma_0$ is its subgraph satisfying
conditions in Theorem \ref{mth1'} so the parabolic ideal
$I_{\Gamma_0}$ is $H_c$-invariant. We will assume that there are at
least two vertices in $\Gamma \setminus \Gamma_0$ so that the
restricted CM system \mref{CMres} is at least two-dimensional hence
non-trivial. Then in all the cases below the parabolic strata are in
one to one correspondence with the isomorphism classes of the
Coxeter subgraphs $\Gamma_0$ except the following cases for the
$\mathcal E_7$ root system (\cite{OS}). Namely, there are two
different $\mathcal A_1^3$ strata denoted as $(\mathcal E_7,
\mathcal A_1^3)_{1,2}$ and two different $\mathcal A_5$ strata
denoted as $(\mathcal E_7, \mathcal A_5)_{1,2}$.

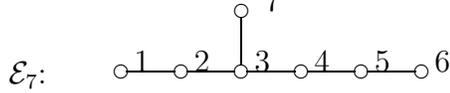
\setlength\unitlength{0.8cm}

\newcommand{\Esevena}{
\begin{picture}(8,1.5)(0,1.8)
\put(1,2){\circle{0.2}} \put(1.1,2){\line(1,0){.8}}
\put(2,2){\circle*{0.2}} \put(2.1,2){\line(1,0){.8}}
\put(3,2){\circle{0.2}} \put(3.1,2){\line(1,0){.8}}
\put(4,2){\circle*{0.2}} \put(4.1,2){\line(1,0){.8}} \put(5
,2){\circle{0.2}} \put(3  ,2.1){\line(0,1){.8}} \put(3
,3){\circle{0.2}} \put(5.1,2){\line(1,0){.8}}
\put(6,2){\circle*{0.2}}
\end{picture}
}
\newcommand{\Esevenb}{
\begin{picture}(8,1.5)(0,1.8)
\put(1,2){\circle{0.2}1}\put(1.1,2){\line(1,0){.8}}
\put(2,2){\circle{0.2}2} \put(2.1,2){\line(1,0){.8}}
\put(3,2){\circle{0.2}3} \put(3.1,2){\line(1,0){.8}}
\put(4,2){\circle{0.2}4} \put(4.1,2){\line(1,0){.8}} \put(5
,2){\circle{0.2}5} \put(3  ,2.1){\line(0,1){.8}} \put(3
,3){\circle{0.2} 7}  \put(5.1,2){\line(1,0){.8}}
\put(6,2){\circle{0.2}6} 
\end{picture}
}

$\mathcal E_7$: $\Esevenb$ 

The stratum $(\mathcal E_7, \mathcal A_1^3)_2$ corresponds to the
choice of subgraph $\Gamma_0$ corresponding to the roots numbered by
4,6,7 in the diagram. The stratum $(\mathcal E_7, \mathcal A_1^3)_2$
corresponds to any other choice of the subgraph of type $\mathcal
A_1^3$. Also the stratum $(\mathcal E_7, \mathcal A_5)_2$
corresponds to the choice of the subgraph with the vertices
3,4,5,6,7. And the stratum  $(\mathcal E_7, \mathcal A_5)_1$
corresponds to any of the two remaining embeddings of the subgraph
$\mathcal A_5$ into $\mathcal E_7$.

In the following table for every system $(\Gamma,\Gamma_0)$ we
specify vectors $A$ in this system, their multiplicities, dimension
of the linear space spanned by $A$, and the number of vectors in $A$.

\vspace{10mm}
\newpage

\begin{longtable}{|c|c|p{83mm}|p{13mm}|c|p{6mm}|}

\hline   & $(\Gamma,\Gamma_0)$   &  Vectors $ A$ of the restricted
CM system  & Multipl. &
Dim & $| A|$  \hspace{2cm}  \\

\hline 1 &  $(\mathcal E_8, \mathcal A_1)$ & $\pm e_1\pm e_2 \pm e_3
\pm e_4 \pm e_5 \pm e_6 + \sqrt{2}e_7$ (even $\natural$ of minuses);
$e_i \pm e_j (1\le i<j\le 6)$; $e_7$
&1/2 &7&91 \\
&& $e_1\pm e_2 \pm e_3 \pm e_4 \pm e_5 \pm e_6$ (odd $\natural$ of
minuses); $\sqrt{2}(e_i \pm e_7) (1\le i\le 6)$
 &1& & \\

\hline 2 & $({\mathcal E}_8, {\mathcal A}_1^2)$ & $e_1\pm e_2 \pm e_3 \pm e_4 \pm e_5 \pm
e_6$ & 1
   & 6 & 68 \\
&&$e_1, e_2, e_3, e_4, e_5, e_6$ & 2 &&\\
&&$ e_i \pm e_j (1\le i<j\le 6)$ &1/2 && \\

\hline 3 & $(\mathcal E_8, \mathcal A_2)$ & $e_1\pm e_2 \pm e_3 \pm
e_4 \pm e_5 \pm \sqrt{3} e_6$ (even $\natural$ of minuses);
 $ e_i \pm e_j (1\le i<j\le 5)$ &1/3  & 6  & 63 \\
&& $e_1\pm e_2 \pm e_3 \pm e_4 \pm e_5 \pm \frac1{\sqrt{3}}e_6$ (odd
$\natural$ of minuses); $\sqrt{3}e_i \pm e_6 (1\le i \le 5)$; $e_6$
& 1
&&\\

\hline 4 & $(\mathcal E_8,\mathcal A_1^3)$ & $e_1\pm e_2 \pm e_3 \pm
e_4 \pm \sqrt{2}
e_5$; $\sqrt{2}e_i \pm e_5 (1\le i \le 4)$ & 1 & 5 & 49  \\
&& $ e_1\pm e_2 \pm e_3 \pm e_4$; $e_1, e_2, e_3, e_4$ &2&& \\
&& $ e_i \pm e_j \,\, (1\le i<j\le 4)$ & 1/2 &&\\
&& $ e_5$ & $9/2$&&  \\

\hline 5 & $(\mathcal E_8,\mathcal A_3)$ & $ e_1\pm e_2 \pm e_3 \pm e_4 \pm e_5$ &1 & 5 & 41 \\
&&  $e_i \pm e_j \,  (1\le i<j\le 5)$ & 1/4 && \\
&& $e_1, e_2, e_3, e_4, e_5$ & 3/2 &&\\

\hline 6 & $(\mathcal E_8,\mathcal A_1^4)$ & $e_1\pm e_2 \pm  e_3$;
$e_1\pm e_2 \pm 
e_4$; $e_1\pm e_3 \pm  e_4$; $e_2\pm e_3 \pm  e_4$ & 1&  4 &  32\\
&& $e_1\pm e_2$; $e_1\pm e_3$; $e_1\pm e_4$; $e_2\pm e_3$; $e_2\pm
e_4$; $e_3\pm e_4$ &2&&\\
&& $e_1$, $e_2$, $e_3$, $e_4$ & 9/2 &&\\

\hline 7 & $(\mathcal E_8,\mathcal A_2^2)$ & $\pm e_1 \pm e_2 +
\frac{1}{\sqrt{3}} e_3 \pm \sqrt{3} e_4$ (odd $\natural$ of
minuses); $\pm e_1 \pm e_2 \pm \sqrt{3} e_3 + \frac{1}{\sqrt{3}}
e_4$ (odd $\natural$ of minuses); $\sqrt{3}e_1 \pm e_3; \sqrt{3} e_2
\pm e_3; \sqrt{3} e_1 \pm e_4; \sqrt{3} e_2 \pm e_4; e_3,
e_4$ & 1 & 4 & 30 \\
&& $e_1\pm e_2
\pm \sqrt{3} e_3 \pm \sqrt{3} e_4$ (even $\natural$ of minuses); $e_1 \pm e_2$ & 1/3 && \\
&& $e_1+e_2 \pm \frac{1}{\sqrt{3}}(e_3+e_4)$; $e_1-e_2 \pm
\frac{1}{\sqrt{3}}(e_3-e_4)$; $e_3 \pm e_4$ & 3 && \\

\hline 8 & $(\mathcal E_8,\mathcal A_4)$ & $e_1\pm e_2 \pm e_3 \pm
\sqrt{5} e_4$ (even $\natural$ of minuses); $e_1 \pm e_2$; $e_1 \pm
e_3$; $e_2 \pm e_3$
& 1/5 & 4 & 25 \\
&& $e_1\pm e_2 \pm e_3 \pm \frac{1}{\sqrt{5}} e_4$ (even $\natural$
of minuses); $e_5$ & 2 &&\\
&& $e_1\pm e_2 \pm e_3 \pm \frac{3}{\sqrt{5}} e_4$ (odd $\natural$
of minuses); $\sqrt{5} e_1 \pm e_4; \sqrt{5} e_2 \pm e_4$; $\sqrt{5}
e_3 \pm e_4$ &1&& \\

 \hline 9 & $\mathcal (\mathcal E_8,\mathcal D_4)$ & $\mathcal F_4$ & $4/3, 1/6$ & 4 & 24\\
\hline

\end{longtable}

\newpage

\begin{longtable}{|c|c|p{83mm}|p{13mm}|c|p{6mm}|}

\hline   & $(\Gamma,\Gamma_0)$   &  Vectors $ A$ of the restricted
CM system  & Multipl. &
Dim & $| A|$  \hspace{2cm}  \\

\hline 10 & $(\mathcal E_8,\mathcal A_5)$ & $ e_1 + e_2 \pm \sqrt{6}
e_3$; $e_1+e_2$
&1/6&  3 & 13 \\
&& $e_1-e_2 \pm \frac{2\sqrt{6}}{3}e_3$; $\sqrt{6} e_1 \pm e_3$;
$\sqrt{6} e_2 \pm e_3$ & 1 && \\
&& $e_1 + e_2 \pm \frac{\sqrt{6}}{3}e_3$; $e_3$ & 5/2 && \\
&& $e_1-e_2$ & 7/2&&\\

\hline 11 & $(\mathcal E_8,\mathcal D_5)$ & $e_1 \pm e_2 \pm e_3$  &2  & 3 & 13 \\
&&  $e_1 \pm e_2$; $e_1 \pm e_3$; $e_2 \pm e_3$ & 1/8 &&  \\
&& $e_1, e_2, e_3$  & 5/4 && \\

\hline 12 & $(\mathcal E_8,\mathcal A_3^2)$ & $e_1, e_2$ & 15/2 & 2 & 8 \\
&& $e_1 \pm 2 e_2$; $e_1 \pm \frac12 e_2$ & 1 && \\
&& $e_1 \pm e_2$ & 4 && \\

\hline 13 & $(\mathcal E_8,\mathcal A_6)$ & $e_1$ & 1/7 & 2 & 6 \\
&& $e_2$ & 6 && \\
&& $\sqrt{7} e_1 \pm 3 e_2$ & 1 && \\
&& $\sqrt{7} e_1 \pm e_2$ & 3 &&\\

\hline 14 & $(\mathcal E_8,\mathcal D_6)$ & $\mathcal B_2$ & 33/10, 6/5 & 2 & 4 \\

\hline 15 & $(\mathcal E_8,\mathcal E_6)$ & $\mathcal G_2$ & 23/12, 1/12 & 2 & 6 \\

\hline


\hline

\hline 16 & $(\mathcal E_7,\mathcal A_1)$ & $ \sqrt{2}e_6 \pm
\sqrt{2} e_5 \pm e_1 \pm e_2 \pm e_3 \pm e_4$ (odd $\natural$ of
minuses in the last four terms); $ e_i \pm e_j$ $(1\le i<j\le 4)$;
$e_5, e_6$ & 1/2 & 6
& 46 \\
&& $\sqrt{2} e_6 \pm e_1 \pm e_2 \pm e_3 \pm e_4$ (even $\natural$
of minuses);  $\sqrt{2}e_i \pm e_5$ $(1\le i \le 4)$ &1 && \\

\hline 17 & $(\mathcal E_7,\mathcal A_1^2)$ & $\sqrt{2}e_5 \pm e_1
\pm e_2 \pm
e_3 \pm e_4$ &1 & 5 & 33 \\
&& $ e_i \pm e_j$ $(1\le i<j\le 4)$; $e_5$ & 1/2 &&\\
&& $e_1$, $e_2$, $e_3$, $e_4$ &2&&
\\

\hline 18 & $(\mathcal E_7,\mathcal A_2)$ & $ \pm e_1 \pm e_2 \pm
e_3 \pm \sqrt{3} e_4 + \sqrt{2} e_5$ (odd $\natural$ of minuses);
$e_1 \pm e_2$; $e_1 \pm e_3$; $e_2 \pm
e_3$; $e_5$ & 1/3& 5 & 30 \\
&&  $ \pm e_1 \pm e_2 \pm e_3 \pm \frac1{\sqrt{3}} e_4 + \sqrt{2}
e_5$
(even $\natural$ of minuses); $\sqrt{3}e_1 \pm e_4$; $\sqrt{3}e_2 \pm e_4$; $\sqrt{3}e_3 \pm e_4$; $e_4$ & 1&&\\

\hline 19 & $(\mathcal E_7,\mathcal A_1^3)_1$ & $\mathcal F_4$ & 2,
1/2 & 4 & 24
\\
\hline

\end{longtable}

\newpage

\begin{longtable}{|c|c|p{83mm}|p{13mm}|c|p{6mm}|}

\hline   & $(\Gamma,\Gamma_0)$   &  Vectors $ A$ of the restricted
CM system  & Multipl. &
Dim & $| A|$  \hspace{2cm}  \\

\hline 20 & $(\mathcal E_7,\mathcal A_1^3)_2$ & $e_1 \pm e_3 \pm
e_4$;
$e_1 \pm e_2 \pm e_4$; $e_2 \pm e_3 \pm e_4$ & 1& 4 & 22 \\
&& $e_1 \pm e_2$; $e_1 \pm e_3$; $e_2 \pm e_3$ &2 && \\
&& $e_1$, $e_2$, $e_3$ & 1/2 &&\\
&& $e_4$ & 9/2 &&\\
\hline

\hline 21 & $(\mathcal E_7,\mathcal A_3)$ & $
e_1 \pm e_2 \pm e_3 \pm \sqrt{2} e_4$ & 1& 4 & 18 \\
&&  $ e_1\pm e_2$; $ e_1\pm e_3$; $ e_2\pm e_3$; $e_4$ & 1/4 &&\\
&& $e_1$, $e_2$, $e_3$ & 3/2&&\\

 \hline 22 & $(\mathcal E_7,\mathcal A_1^4)$ & $e_1 \pm e_2 \pm  e_3$ &1 &  3 & 13 \\
&&   $e_1\pm e_2$; $ e_1\pm e_3$; $ e_2\pm e_3$ &2&&\\
&& $e_1, e_2, e_3$ & 9/2 &&\\

\hline 23 & $(\mathcal E_7,\mathcal A_2^2)$ & $\sqrt{3}(e_1-e_2) \pm
\sqrt{2}e_3 $;
$e_3$ &1/3 &  3 & 13 \\
&& $e_1 + 3 e_2 \pm \sqrt{6}e_3$; $3
e_1 +  e_2 \pm \sqrt{6}e_3$; $e_1$, $e_2$ & 1&&\\
&& $e_1-e_2 \pm \sqrt{6} e_3$; $e_1 \pm e_2$ & 3 &&\\

\hline 24 & $(\mathcal E_7,\mathcal A_4)$ & $e_1 - \sqrt{5} e_2 \pm
\sqrt{2}e_3$;
$e_3$ & 1/5 & 3 & 10\\
&& $e_1 +\frac{3}{\sqrt{5}} e_2 \pm \sqrt{2}e_3$; $\sqrt{5}e_1 \pm
e_2$ & 1 && \\
&& $e_1 - \frac1{\sqrt{5}} e_2 \pm \sqrt{2}e_3$; $e_2$ & 2 &&\\

\hline 25 & $(\mathcal E_7,\mathcal D_4)$ & $\mathcal B_3$ & $4/3, 1/6$ & 3&9\\

\hline 26 & $(\mathcal E_7,\mathcal D_5)$ & $e_1 \pm \sqrt{2}e_2$ & 2 & 2&4\\
&& $e_1$ & 5/4 &&\\
&& $e_2$ & 1/8 &&\\

\hline 27 & $(\mathcal E_7,\mathcal A_5)_1$ & $ \mathcal C_2(\frac52, \frac72)$ & 1, 5/2, 7/2 & 2&4\\

\hline 28 & $(\mathcal E_7,\mathcal A_5)_2$ & $\mathcal G_2$ & $5/2, 1/6$ & 2&6\\

\hline

\hline 29 & $(\mathcal E_6,\mathcal A_1)$ & $\pm e_1 \pm e_2 \pm e_3
+ \sqrt{2} e_4 \pm \sqrt{3} e_5$ (even $\natural$ of minuses); $e_1
\pm e_2$; $e_1
\pm e_3$; $e_2 \pm e_3$; $e_4$ & 1/2 &  5 &  25 \\
&& $\pm e_1 \pm e_2 \pm e_3 + \sqrt{3} e_5$ (odd $\natural$ of
minuses); $\sqrt{2}e_1 \pm e_4$; $\sqrt{2}e_2 \pm e_4$; $\sqrt{2}e_3
\pm e_4$ & 1 &&\\

\hline 30 & $(\mathcal E_6,\mathcal A_1^2)$ & $ e_1 \pm e_2 \pm e_3 \pm \sqrt{3} e_4$ & 1 & 4 &  17 \\
&& $e_1 \pm e_2$; $e_1 \pm e_3$; $e_2 \pm e_3$ & 1/2 && \\
&& $e_1, e_2, e_3$ & 2 &&\\
\hline

\end{longtable}

\newpage

\begin{longtable}{|c|c|p{83mm}|p{13mm}|c|p{6mm}|}

\hline   & $(\Gamma,\Gamma_0)$   &  Vectors $ A$ of the restricted
CM system  & Multipl. &
Dim & $|A|$  \hspace{2cm}  \\

\hline 31 & $(\mathcal E_6,\mathcal A_2)$ & $ \pm e_1 \pm e_2 \pm
\sqrt{3} e_3 + \sqrt{3} e_4$ (even  $\natural$ of minuses); $e_1 \pm
e_2$ & 1/3 & 4
& 15 \\
&& $ \pm e_1 \pm e_2 \pm \frac1{\sqrt{3}} e_3 + \sqrt{3} e_4$ (odd
$\natural$ of minuses); $\sqrt{3} e_1 \pm e_3$; $\sqrt{3} e_2 \pm
e_3$; $e_3$ &1 && \\

\hline 32 & $(\mathcal E_6,\mathcal A_1^3)$ & $e_1 \pm \sqrt{2} e_2
\pm \sqrt{3}e_3$;
$\sqrt{2} e_1 \pm e_2$ & 1 & 3 & 10 \\
&& $e_1 \pm \sqrt{3} e_3$; $e_1$ & 2 && \\
&& $e_2$ & 9/2 && \\

\hline 33 & $(\mathcal E_6,\mathcal A_3)$ & $\mathcal C_3(1/4,0)$ & 1, 1/4,  3/2 & 3& 8 \\


\hline 34 & $(\mathcal E_6,\mathcal A_2^2)$ & $\mathcal G_2$ & 3, 1/3 &  2 &  6 \\

\hline 35 & $(\mathcal E_6, \mathcal D_4)$ & $\mathcal A_2$ &4/3 &  2 & 3 \\

\hline 36 & $(\mathcal E_6,\mathcal  A_4)$ & $ e_1$ & 1/5 &  2 & 4 \\
&& $e_2$ & 1 && \\
&& $\sqrt{5}e_1 \pm \sqrt{3} e_2$ & 2 && \\


\hline

\hline 37 &  $(\mathcal H_4,\mathcal A_1)$ & $(\sqrt{5}+1)e_1 \pm
2e_2 \pm (\sqrt{5}-1)e_3$; $ 2e_1 \pm (\sqrt{5}-1)e_2 \pm
(\sqrt{5}+1)e_3$; $(\sqrt{5}-1)e_1  \pm (\sqrt{5}+1)e_2 \pm 2e_3$;
$e_1$, $e_2$, $e_3$
  &1/2&  3 &  31 \\
&& $2e_1  \pm (\sqrt{5}+3)e_3$; $(\sqrt{5}+3)e_2  \pm 2e_3$;
$(\sqrt{5}+3)e_1  \pm 2e_2$; $e_1 \pm e_2 \pm e_3$ & 1&&\\
&& $2e_1 \pm (\sqrt{5}+1)e_2$; $(\sqrt{5}+1)e_1 \pm 2 e_3$; $2e_2
\pm (\sqrt{5}+1)e_3$ &2&&\\

\hline 38 &  $(\mathcal H_4,\mathcal A_2)$ & $(\sqrt{5}\pm2)e_1+\sqrt{3}e_2$; $\sqrt{5}e_1-\sqrt{3}e_2$ & 1/3 &  2& 12\\
&& $\sqrt{3}e_1+\sqrt{5}e_2$, $\sqrt{3}e_1 -(\sqrt{5}\pm2)e_2$
&4&&\\
&& $\sqrt{3}e_1+(\sqrt{5}\pm 4)e_2$; $\sqrt{3}e_1\pm e_2$,
$\sqrt{15}e_1+e_2$, $e_2$ &1&&\\

\hline 39 &  $(\mathcal H_4,\mathcal I_2(5))$ & ${\mathcal I}_2(10)$ & 1/5, 2 & 2 & 10\\

\hline 40 &  $(\mathcal H_4,\mathcal A_1^2)$ & $e_1, e_2$ & 13/2 & 2 & 12\\
&& $e_1 \pm e_2$; $(\sqrt{5}\pm 1)e_1 \pm 2e_2$ &2&&\\
&& $(\sqrt{5}\pm 3)e_1 \pm 2 e_2$ & 1 &&\\

 \hline

\hline 41 &  $(\mathcal H_3,\mathcal A_1)$ & $e_1$, $e_2$ & 1/2 &  2& 6\\
&& $(\sqrt{5}+1) e_1 \pm 2 e_2$ &2&&\\
&& $(\sqrt{5}+3) e_1 \pm 2 e_2$ &1&&\\

\hline
\end{longtable}

The restrictions of Coxeter root systems appear also, in particular,
in the context of $\vee$-systems \cite{FV}. We note that the number
of vectors in the $\vee$-system $(\mathcal H_4, \mathcal A_1)$ is
stated inaccurately in \cite{FV} as some of the vectors listed there
are actually proportional. We also refer to \cite{Goodwin} where, in
particular, bases of the restricted Weyl root systems are discussed.

The CM systems corresponding to the Coxeter restrictions $(\mathcal
E_7, \mathcal D_5), (\mathcal E_6, \mathcal A_4)$ belong to the
one-parametric family of two-dimensional integrable CM systems
introduced and studied in \cite{Taniguchi2}. The system 12 belongs
to the family from Proposition \ref{propos} and the system 22
belongs to the family from Theorem \ref{theoremf4a1}.

We note that in the table above in all the cases there are
non-integer multiplicities. So all integrable systems with integer
multiplicities which are restrictions of the CM systems with Coxeter
root systems already appeared in \cite{CFV3}. Also according to the
result from \cite{Taniguchi1} the group generated by reflections
along the hyperplanes with non-integer multiplicities is finite.
This agrees with the table above.


\section{Invariant parabolic ideals for the complex reflection groups}
\label{sec6}

Let $W$ be an irreducible finite complex reflection group. Let $V$ be its reflection representation. Let $\cal A$ be the set of reflection hyperplanes. For any hyperplane $H \in {\cal A}\subset V$ let $m_H$ be the order of the stabilizer of $H$ in the group $W$, and let $\a_H$ be a covector vanishing on $H$. Let $s_{H,i}$, $i=1,\ldots,m_H-1$ be the set of reflections in $W$ which fix $H$. We numerate these reflections so that $s_{H,i}=s_{H,1}^i$, and we suppose that $det s_{H,1} = \xi_H = e^{2\pi i/m_H}$. Put $s_{H,0}=e \in W$. Let $a(s_{H,i})=a_{H,i}$ be a $W$-invariant function on the set of reflections.

For any reflection $s_{H,i}$ we choose a pair $\a_{H,i} \in V^*$, $\a_{H,i}^\vee \in V$ such that $s_{H,i}(f) = f - f(\a_{H,i}^\vee)\a_{H,i}$ for any $f \in V^*$. Note that these pairs are not uniquely defined by reflections but the elements $\a_{H,i}\otimes \a_{H,i}^\vee \in V^*\otimes V$ are. Define the bilinear form $B:V^* \otimes V \to \C$ by the formula
$$
B(f,v)= \sum_{H \in \cal A} \sum_{i=1}^{m_H-1} a_{H,i} \alpha_{H,i}(v) f(\a_{H,i}^\vee).
$$
Because of $W$-invariance one has $B(f,v)= h_{W, a} f(v)$ for some constant $h_{W,a}$. More directly this coefficient is defined as
\beq{gencoxc}
 h_{W,a} = \sum_{H \in \cal A} \sum_{i=1}^{m_H-1} a_{H,i} (1-s_{H,i}),
\eeq
where the reflections $s_{H,i}$ act in $V^*$ or, equivalently, in $V$.
\begin{remark}
It is noted in \cite{GG} that $h_W=\sum_{H\in\cal A} \sum_{i=1}^{m_H-1} (1- s_{H,i}) = \frac{1}{\dim V} \sum_{H \in \cal A} m_H$, so in particular $h_W$ is a generalization of the Coxeter number for the complex reflection group $W$. In the case $a_{H,i}=const = a$ we have $h_{W,a}= a h_W$.
\end{remark}
Recall that the Dunkl operators are defined for any $\xi \in V$ as \cite{DO}
\beq{dunklcomplex}
\nabla_\xi = \partial_\xi - \sum_{H \in \cal A} \frac{\a_H(\xi)}{\a_H} \sum_{i=1}^{m_H-1}  a_{H,i} (1-s_{H,i}).
\eeq
More exactly the definition of the Dunkl operator in \cite{DO} is
\beq{dunklcomplexoriginal}
\nabla_\xi =\partial_\xi- \sum_{H \in \cal A} \frac{\a_H(\xi)}{\a_H} \sum_{t=1}^{m_H-1}  b_{H,t} \sum_{i=0}^{m_H-1} \xi_H^{ i t}s_{H,i},
\eeq
and the formulas \mref{dunklcomplex}, \mref{dunklcomplexoriginal} coincide if the parameters are related by
$$
a_{H,s} = - \sum_{i=1}^{m_H-1} b_{H,i} \xi^{s i}, \,\, s=1,\ldots, m_H-1
$$
for any $H \in \cal A$.

Let now $W_0$ be a parabolic subgroup of $W$ that is $W_0$ is the stabilizer of an intersection  $L$ of the reflection hyperplanes. The corresponding {\it parabolic stratum} is the orbit
$$
D_{W_0}=\bigcup_{w \in W} w(L).
$$
The associated {\it parabolic ideal} $I_{W_0}$ is defined
as the set of polynomials vanishing on the stratum, that is
$I_{W_0}=\{p\in \C[x]\mid p|_{D_{W_0}}=0\}$.
\begin{remark}
In the case of a real reflection group (Section \ref{section-real-parabolic}) the ideal $I_{\Gamma_0}$ coincides with the ideal $I_{W_0}$ where $W_0$ is the parabolic subgroup generated by simple reflections corresponding to the vertices of $\Gamma_0$.
\end{remark}
We are going to determine the parabolic strata
$D_{W_0}$ such that ideals $I_{W_0}$ are invariant under
the rational Cherednik algebra associated to $W$. Equivalently, the ideals are invariant under the Dunkl operators \mref{dunklcomplex}.

Let $W_0=W_1\times\ldots\times W_k$ be the decomposition of $W_0$ into the irreducible parabolic subgroups so that $V=\oplus_{i=1}^k V_i \oplus L$ where each $V_i$ ($1 \le i \le k$) is the reflection representation for $W_i$ and $W_i$ acts trivially in $V_j$, $j\ne i$ and in $L$.

\begin{theorem}\label{mth1-complex}
The parabolic ideal $I_{W_0}$ is invariant under the
Dunkl operators \mref{dunklcomplex} for any $\xi \in V$  if and only if
$h_{W_i,a}=1$ for any $i=1,\ldots,k$.
\end{theorem}
The proof is similar to the proof of Theorem \ref{mth1}. One can check that invariance of the ideal under $\nabla_\xi$ is equivalent to the property that for any $\g \in V^*$ such that $\gamma|_L=0$ one has
$$
\gamma(\xi) =  \sum_{H \in \cal B} \sum_{i=1}^{m_H-1} a_{H,i}\a_{H,i}(\xi) \gamma(\a_{H,i}^\vee),
$$
where $\cal B \subset \cal A$ is the collection of the reflection hyperplanes containing the subspace $L$. Then Thorem \ref{mth1-complex} follows.

\begin{remark} Similar to the real case (Theorem \ref{radidth}) any $H_c(W)$-invariant radical ideal corresponds to the union of parabolic strata so that the parabolic subgroups $W_i$ defining each stratum satisfy the property $h_{W_i,a}=1$.
\end{remark}

Consider now the complex reflection group $W=G(m,p,N)$ and its natural
action in $\C^N$. Recall that the group $G(m,p,N)$ defined when
$p|m$ is generated by the elements $s_{ij}^k$ for $1\le i <j \le N$,
$k=0,\ldots,m-1$ and the elements $\tau_i$ for $i=1,\ldots,N$. The
element $\tau_i$ acts on the basis vectors as $\tau_i(e_i)=\eta^{-1}
e_i$, where $\eta=e^{2\pi i p/m}$ 
and $\tau_i(e_j)=e_j$ for $j\ne i$. The elements $s_{ij}^k$ defined
for $i \ne j$ act as $s_{ij}^k (e_j)=\xi^{k} e_i$, $s_{ij}^k
(e_i)=\xi^{-k} e_j$, where $\xi=e^{2\pi i/m}$, and $s_{ij}^k
(e_s)=e_s$ for $s\ne i,j$. The complex reflections $s_{ij}^k$ are reflections of order 2 at the
hyperplanes $x_i-\xi^k x_j=0$. The complex reflections $\tau_i$ are
reflections of order $m/p$ at the hyperplanes $x_i=0$. We are going to specify invariant parabolic ideals explicitly.

 The Dunkl operators for the complex reflection group $G(m,p,N)$
depend on $m/p$ complex parameters $c_0,\ldots, c_{\frac{m}{p}-1}$
and have the form \cite{DO}
\beq{cdunkl} \nabla_i =
\partial_i-c_0 \sum_{\nad{j=1}{j\ne i}}^N \sum_{k=0}^{m-1}
\frac{1-s_{ij}^k}{x_i-\xi^k x_j}-\sum_{t=1}^{\frac{m}{p}-1} c_t
\sum_{s=0}^{\frac{m}{p}-1}\frac{\eta^{-st}\tau_i^{s}}{x_i}, \eeq
$i=1,\ldots,N$. The commutativity $[\nabla_i,\nabla_j]=0$ holds.

The parabolic strata are
the $G(m,p,N)$-orbits of the intersection of the reflection
hyperplanes $x_i-\xi^k x_j=0$, $x_s=0$.
Consider the plane $\pi_{q,r}^\ve \subset \C^N$ given by the
equations
\begin{multline}\label{pimkbC}
x_1=x_2=\ldots=x_{r},\quad x_{r+1}= x_{r+2}=
\ldots= x_{2 r},\\
\ldots\\
\ve x_{(q-1)r+1}= x_{(q-1)r+2}= \ldots = x_{qr},
\end{multline}
where $\ve^m=1$. We may assume that $\ve=1$ unless $qr = N$. The
corresponding parabolic stratum $D_{q,r}^\ve\subset \C^N$ is the
orbit 
\beq{stratumBN1C} D_{q,r}^\ve=\bigcup_{w\in { G}(m,p,N)}
w(\pi_{q,r}^\ve).\eeq
\begin{prop}\label{BN1C}
The parabolic ideal
$$I_{q,r}^\ve=\{f\in \C[x_1,\ldots,x_N] \ | \ f|_{D_{q,r}^\ve}=0\}$$
corresponding to the stratum \mref{pimkbC}, \mref{stratumBN1C} is
invariant under the $G(m,p,N)$ Dunkl operators \mref{cdunkl} if and
only if $c_0=1/r$.
\end{prop}
{\bf Proof.}
The stratum \mref{stratumBN1C} corresponds to the parabolic subgroup of type $W_0 \cong A_{r-1}^q$. The Coxeter number of any irreducible component of $W_0$ equals $r$. The statement follows from Theorem \ref{mth1-complex}.

Consider now the plane $\pi_{l}\subset \C^N$ given by the equations
$x_1=x_2=\ldots=x_{l}=0$.  Consider the corresponding parabolic
stratum $D_{l}$ which is the orbit of the plane $\pi_{l}$ under the
group $G(m,p,N)$: \beq{stratumBN2C} D_{l}=\bigcup_{i_1<\ldots <i_l}
\{x| x_{i_1}=\ldots=x_{i_l}=0\}. \eeq Define the ideal $I_{l}=\{f\in
\C[x_1,\ldots,x_N]  \mid f|_{D_{l}}=0\}$.

\begin{prop}\label{BN2C}
The Dunkl operators \mref{cdunkl} for the group $G(m,m,N)$ preserve
the ideal $I_l$ if and only if $c_0=\frac{1}{m(l-1)}$. The Dunkl
operators \mref{cdunkl} for the group $G(m,p,N)$ with $p<m$ preserve
the ideal $I_l$ if and only if $(l-1)c_0+ c_1 p^{-1} =m^{-1}$.
\end{prop}
{\bf Proof.} Let $f \in I_l$. We analyze first the condition
$\nabla_1 f \in I_l$. Since polynomial $f$ vanishes on the plane
$\pi_l: x_1=x_2=\ldots=x_{l}=0$ we can represent $f$ as
$f=\sum_{i=1}^l x_i f_i$ for some polynomials $f_i$. In order to
compute $\nabla_1 f$ we note at first that for $2 \le j \le l$
$$
\left(\frac{1-s_{1j}^k}{x_1-\xi^k x_j} \sum_{i=1}^l x_i f_i\right)
|_{\pi_l}= f_1|_{\pi_l} - \xi^{-k} f_j|_{\pi_l},
$$
and that $\sum_{k=0}^{m-1} \xi^{-k}=0$,  $\sum_{s=0}^{\frac{m}{p}-1}
\eta^{-s t}=0$ for $\frac{m}{p}-1 \ge t \ge 1$.
 Then it follows that
$$
\nabla_1 f|_{\pi_l}=(1-c_0 m (l-1) - c_1 \frac{m}{p})f_1|_{\pi_l},
$$
where we assume that $c_1=0$ for the group $G(m,m,N)$. Therefore
$\nabla_1 f|_{\pi_l}=0$ iff $m^{-1} =c_0 (l-1) + {c_1}{p}^{-1}$.
The property that $\nabla_1 f$ vanishes on other subspaces $\pi$ from
$D_l$ is either satisfied for all values of parameters (when $e_1
\in \pi$) or is satisfied under the same relation among the parameters
as for the plane $\pi_l$ (when the vector $e_1$ is orthogonal to the
subspace $\pi$). This shows that invariance of the ideal $I_l$ under
the Dunkl operator $\nabla_1$, as well as under all operators
$\nabla_i$, is equivalent to the condition $m^{-1} =c_0 (l-1) +
{c_1}{p}^{-1}$ as stated.
\begin{remark}
It follows from Proposition \ref{BN2C} and Theorem \ref{mth1-complex} that in the case $W=G(m,p,N)$ the generalized Coxeter number \mref{gencoxc} can be rearranged as
$$
h_{W,a}=m(N-1)c_0 + \frac{m}{p} c_1,
$$
where $c_0=a(s_{ij}^k)$ and $c_1=\frac{p}{m} \sum_{k=1}^{\frac{m}{p}-1} a(\tau_i^k)(1-\eta^{-k})$, and we assume $c_1=0$ in the case $m=p$.
\end{remark}

Propositions \ref{BN1C}, \ref{BN2C} and their proofs imply the
following
\begin{theorem}\label{theoremlast}
All possible parabolic strata defining parabolic ideals invariant
under the $G(m,p,N)$ Dunkl operators \mref{cdunkl} are either the
strata $D_{q,r}^\ve$ defined by \mref{stratumBN1} with $c_0=1/r$ or
strata $D_l$ defined by \mref{stratumBN2} with parameters values
specified in Proposition \ref{BN2C} or the strata $D_{q,r,l}$. The
latter stratum is the $G(m,p,N)$-orbit of the subspace
\begin{multline}
x_1=x_2=\ldots=x_{r},\quad x_{r+1}= x_{r+2}=
\ldots= x_{2 r},\\
\ldots\\
x_{(q-1)r+1}= x_{(q-1)r+2}= \ldots = x_{qr},
\\ x_{qr+1}=\ldots=x_{qr+l}=0,
\end{multline}
where $r>1$, $l \ge 1$, $qr+l \le N$. The corresponding parameters
satisfy $c_0=1/r$, $(l-1)c_0+ \frac{c_1}{p} =1/m$, where it is
assumed that $c_1=0$ in the case of the group $G(m,m,N)$.
\end{theorem}

To conclude the consideration of the invariant parabolic ideals for
the complex reflection group $G(m,p,N)$ we note that when $N=2$ and
$p$ is even the roots $e_1-\xi^s e_2$, $s=0,\ldots,m-1$ form two
$G(m,p,N)$-orbits. This adds extra parameter $\tilde c_0$  to the
associated Dunkl operators comparing to the general case
\mref{cdunkl}. Namely the operators take the form $$ \nabla_1 =
\partial_1-c_0 \sum_{k=0}^{\frac{m}{2}-1}
\frac{1-s_{12}^{2k}}{x_1-\xi^{2k} x_2}-\tilde c_0
\sum_{k=0}^{\frac{m}{2}-1} \frac{1-s_{12}^{2k+1}}{x_1-\xi^{2k+1}
x_2}-\sum_{t=1}^{\frac{m}{p}-1} c_t
\sum_{s=0}^{\frac{m}{p}-1}\frac{\eta^{-st}\tau_1^{s}}{x_1},$$
\beq{cdunkl2}\nabla_2 =
\partial_2-c_0 \sum_{k=0}^{\frac{m}{2}-1}
\frac{1-s_{12}^{2k}}{x_2-\xi^{2k} x_1}-\tilde c_0
\sum_{k=0}^{\frac{m}{2}-1} \frac{1-s_{12}^{2k+1}}{x_2-\xi^{2k+1}
x_1}-\sum_{t=1}^{\frac{m}{p}-1} c_t
\sum_{s=0}^{\frac{m}{p}-1}\frac{\eta^{-st}\tau_1^{s}}{x_2}.\eeq
Then
we have the following four parabolic ideals
$$I_1=\{f\in \C[x_1,x_2]| f=0 \, \mbox{ if } \,  x_1=\xi^{2k}x_2 \quad \forall k=0,\ldots, \frac{m}{2}-1
\},$$
$$I_2=\{f\in \C[x_1,x_2]| f=0 \, \mbox{ if } \,  x_1=\xi^{2k+1}x_2 \quad \forall k=0,\ldots, \frac{m}{2}-1
\},$$ $I_3=\{f\in \C[x_1,x_2]| f=0 \, \mbox{ if } \,  x_1 x_2=0\}$,
$I_4=\{f\in \C[x_1,x_2]| f=0 \, \mbox{ if } \,  x_1=x_2=0 \}$.

The following Proposition can be established similarly to the
previous results of this Section.

\begin{prop}\label{proplast}
The parabolic ideal $I_1$ is invariant under the Dunkl operators
\mref{cdunkl2} iff $c_0=1/2$, the parabolic ideal $I_2$ is invariant
under the Dunkl operators \mref{cdunkl2} iff $\tilde  c_0=1/2$, the
parabolic ideal $I_3$ is invariant under the Dunkl operators
\mref{cdunkl2} iff $c_1=p/m$, the parabolic ideal $I_4$ is invariant
under the Dunkl operators \mref{cdunkl2} iff $\frac12(c_0+\tilde
c_0)+c_1/p=1/m$.
\end{prop}

\section*{Concluding Remarks}

In the paper we were systematically  deriving generalized CM systems
from the special subrepresentations in the polynomial representation
of the rational Cherednik algebra. The natural development is
to extend this approach to obtain integrable generalizations of the
Calogero-Moser-Sutherland and Macdonald-Ruijsenaars systems
 starting from less degenerate
Cherednik algebras as well as to obtain elliptic generalized Calogero-Moser systems. Some integrable generalizations of these systems
are known from \cite{CFV2}, \cite{SV1}, \cite{Chalykh}, \cite{Feig}, \cite{SV3}, \cite{Hod}. We refer to the recent development in \cite{FS} where generalized Macdonald-Ruijsenaars systems were derived from the special submodules of the polynomial representations of double affine Hecke algebras.

It would also be interesting to analyze if the approach can be
extended to cover matrix integrable systems of the generalized
Calogero-Moser type \cite{CGV}.
A close direction is to investigate possible generalizations of
matrix and scalar CM systems associated to special complex reflection groups (see \cite{Young} and \cite{EFMV}  respectively).

Regarding representations of the rational Cherednik algebras, it is
clear by Theorem \ref{mth1} how to form the chains of submodules in
the polynomial module of the rational Cherednik algebras for any
Coxeter group. In contrast to the ${\mathcal A}_N$ case in general there may be
non-isomorphic Coxeter subgroups with equal Coxeter numbers. Therefore there are non-trivial
intersections of the corresponding parabolic ideals. It would be interesting to investigate if the
subsequent quotients  in the natural chains are irreducible like in the ${\mathcal A}_N$ case
(\cite{K1}, \cite{Enomoto}, \cite{ES}). Also in the paper we consider
submodules corresponding to special singular values only, it would
also be interesting to see if the submodules for other singular
values can be described in a natural way (c.f. \cite{ES} and Proposition
\ref{onedimmult} above).

{\bf Acknowledgements.} I am very grateful to Yu. Berest for
stimulating discussions and to A.P. Veselov for the attention to the
work. I am also very grateful to O.A. Chalykh, V. Dotsenko, I. Gordon, S. Griffeth, C.
Korff, A.N. Sergeev, A. Silantyev, J.T. Stafford, I. Strachan  for useful
discussions. I would like to thank a referee for a few suggestions to improve the paper and especially for the hint on how to extend the description of invariant parabolic ideals for $H_c(G(m,p,N))$ to the case of an arbitrary complex reflection group (Theorem~\ref{mth1-complex}).

The work was partially supported by the EPSRC grant
EP/F032889/1, by the ESF Network ENIGMA (contract
MRTN-CT-2004-5652), by the EPSRC Network grant EP/F029381/1, by British Council (PMI2 Research Cooperation Award).

\vspace{5mm}


\vspace{1cm}

\end{document}